\documentclass[reqno]{amsart}
\usepackage{amsmath,amsfonts,amssymb,amsthm,enumitem,graphicx}
\usepackage[usenames,dvipsnames]{xcolor}
\usepackage[width=\linewidth]{caption}
\usepackage{subcaption}
\usepackage[colorlinks=true,linkcolor=blue]{hyperref} 
\newtheorem{theorem}{Theorem}[section]
\newtheorem{lema}[theorem]{Lemma}
\newtheorem{prop}[theorem]{Proposition}
\newtheorem{coro}[theorem]{Corollary}
\theoremstyle{definition}
\newtheorem{definition}[theorem]{Definition}

\newtheorem{remark}[theorem]{Remark}

\newtheorem{notas}[theorem]{Remarks}
\numberwithin{equation}{section}

\newcommand{\R}{\mathbb{R}}

\newcommand{\A}{\mathbb{A}}

\newcommand{\W}{\Omega}
\newcommand{\w}{\omega}

\newcommand{\ep}{\varepsilon}

\newcommand{\wit}{\widetilde}
\newcommand{\n}[1]{\| #1 \|}
\newcommand{\lsm}{\left[\!\begin{smallmatrix}}
\newcommand{\rsm}{\end{smallmatrix}\!\right]}

\newcommand{\des}{\displaystyle}

\DeclareMathOperator{\cls}{cls} 
\DeclareMathOperator{\Int}{Int} 
 
\definecolor{col}{rgb}{0,0,0.6}
\begin{document}
\title[Exponential ordering with applications to Nicholson systems]
{The exponential ordering for non-autonomous delay systems with applications to compartmental Nicholson systems}
\author[S. Novo]{Sylvia Novo}
\author[R. Obaya]{Rafael Obaya}
\author[A.M. Sanz]{Ana M. Sanz}
\author[V.M.~Villarragut]{V\'{\i}ctor M. Villarragut}
\address[S. Novo]{Departamento de Matem\'{a}tica Aplicada, Escuela de Ingenier\'{\i}as Industriales (Sede Doctor Mergelina), Universidad de
  Valladolid, 47011 Valladolid, Spain.}
  \email{sylvia.novo@uva.es}
\address[R. Obaya]{Departamento de Matem\'{a}tica
Aplicada, Escuela de Ingenier\'{\i}as Indus\-tria\-les (Sede Doctor Mergelina), Universidad de Valladolid,
47011 Valladolid, Spain, and member of IMUVA, Instituto de Investigaci\'{o}n en
Matem\'{a}ticas, Universidad de Va\-lla\-dolid, Spain.}
 \email{rafael.obaya@uva.es}
\address[A.M. Sanz]{Departamento de Did\'{a}ctica de las Ciencias Experimentales, Sociales y de la Matem\'{a}tica,
Facultad de Educaci\'{o}n, Universidad de Valladolid, 34004 Palencia, Spain,
and member of IMUVA, Instituto de Investigaci\'{o}n en  Mate\-m\'{a}\-ti\-cas, Universidad de
Valladolid.} \email{anamaria.sanz@uva.es}
\address[V.M.~Villarragut]{Departamento de Matem\'{a}tica Aplicada a
  la Ingenier\'{\i}a Industrial, Universidad Polit\'{e}cnica de Madrid, Calle de
  Jos\'{e} Guti\'{e}rrez Abascal 2, 28006 Madrid, Spain.}
  \email{victor.munoz@upm.es}
\thanks{The first three authors were partly supported by MICIIN/FEDER project
RTI2018-096523-B-I00 and by Universidad de Valladolid under project PIP-TCESC-2020. The fourth author was partly supported by MICINN/FEDER under projects RTI2018-096523-B-I00 and PGC2018-097565-B-I00.
}
\date{}
\begin{abstract}
The exponential ordering is exploited in the context of non-auto\-no\-mous delay
systems, inducing monotone skew-product semiflows under less restrictive
conditions than usual. Some dynamical concepts linked to the order, such as
semiequilibria, are considered for the exponential ordering, with implications
for the determination of the presence of uniform persistence or the existence of
global attractors. Also, some important conclusions on the long-term dynamics
and attraction are obtained for monotone and sublinear delay systems for this
ordering. The results are then applied to almost periodic Nicholson systems and
new conditions are given for the existence of a unique almost periodic positive
solution which asymptotically attracts every other positive solution.
\end{abstract}
\keywords{Non-autonomous dynamical systems, exponential ordering, uniform
persistence, global attractor, almost periodic Nicholson systems} \subjclass{37C60, 37C65, 37C75, 92D25}
\renewcommand{\subjclassname}{\textup{2020} Mathematics Subject Classification}
\maketitle
\section{Introduction}\label{sec-intro}
The modern dynamical theory of monotone skew-product semiflows has been extensively investigated in the last few years. References such as  Chueshov~\cite{book:chue}, Jiang and Zhao~\cite{jizh}, Novo et al.~\cite{noos07JDE}, and Shen and Yi~\cite{shyi}, among others, are important for the initial development of this theory. In this setting, the skew-product semiflow $\tau: \R^+\times\Omega \times X \to \Omega \times X$ is defined on a product bundle space $\Omega \times X$, where $\Omega$ is a compact metric space under the action of a continuous minimal flow $\sigma: \R\times \Omega \to \Omega$  and $X$ is an ordered Banach space, that is, there is a cone  of positive vectors $X_+$ inducing a partial order relation. Monotonicity means that the order is preserved under the semiflow. Frequently, a first condition required in this theory is that the positive cone has a nonempty interior. A second relevant condition is the uniform stability of the relatively compact semitrajectories, which implies their convergence to minimal subsets of the product space. In absence of uniform stability, the omega-limit sets show  ingredients of dynamical complexity, including even chaotic dynamics.\smallskip

This paper provides a contribution to the dynamical theory of monotone skew-product semiflows generated by non-autonomous functional differential equations (FDEs for short) with finite delay. We consider a family of finite-delay FDEs  over $\W$, $y'(t)= F(\omega{\cdot}t, y_t)$, $\w\in\W$,  where $F:\Omega \times C([-r,0],\R^m) \to \R^m$ is continuous and of class $C^1$  with respect to the second component. Here, $X$ is the Banach space $C([-r,0],\R^m)$ with positive cone $X_+$ given by the componentwise nonnegative functions in $X$. In this paper we investigate those FDEs which generate a semiflow $\tau$  monotone for the exponential ordering  defined  by a quasipositive constant matrix $B$. This ordering, $\leq _B$,   introduced  by Smith and Thieme~\cite{smth90,smth91},  has a positive cone $K_B$ with empty interior in $X$. An  option to overcome this problem is to  consider the definition of $\tau$ on $\R^+\times\Omega \times X_L$, taking the state space $X_L \subset  X$ of Lipschitz functions, even though the continuity of $\tau$ fails for times in the interval~$[0,r]$.\smallskip

In the first part of the paper, we show how important ingredients of the theory
of monotone skew-product semiflows can be extended to this context despite the
difficulties mentioned above. Following the ideas and methods of Chueshov~\cite{book:chue} and Novo et al.~\cite{nono2}, sub and super-equilibria
--now for the exponential ordering-- prove to be useful dynamical objects.
Simple conditions that imply the existence of these functions are provided,
which ensures the applicability of the results. We use arguments of the theory
of non-autonomous dissipative dynamical systems, in the terms given in Kloeden
and Rasmussen~\cite{klra}, Carvalho et al.~\cite{book:CLR}, and Cheban et al.~\cite{chks}, to investigate the
long-term behaviour of the trajectories of the semiflow. We also deal with the
property of uniform persistence, which is an important notion in dynamical
systems that appears in the literature through several formulations. In this
paper, it will refer to  persistence for the order, as introduced in Faria and R\"{o}st~\cite{faro} and Novo et al.~\cite{noos13}. When $\tau$ is globally defined on $\W\times K_B$, we provide criteria based on the existence of semiequilibria to  deduce, firstly, the uniform persistence  of the semiflow in $\Int \wit K_B$ and, secondly, the existence of a global attractor inside the open set $\W\times \Int \wit K_B$.\smallskip

Special attention is paid to the case in which the function $F$ is in addition
sublinear for the exponential ordering. Thanks to the monotonicity and
sublinearity properties of $F$, the semiflow $\tau$ is globally defined on
$\W\times K_B$ and it is sublinear for $\leq_B$. Besides, if there exists a positive bounded semitrajectory that is uniformly strongly above $0$ for the exponential ordering, we deduce that the semiflow is uniformly stable on the relatively compact subsets of $\Omega \times K_B$ that are  uniformly strongly above $0$. As a consequence,  the omega-limit set of every $(\omega, \phi)$  with $\phi \gg_B  0$ is a uniformly stable and strongly positive minimal set which admits a fiber distal flow extension. Furthermore, if $E$ denotes the union of all these omega-limit sets, we show that $E \subset \Omega  \times \Int \wit K_B$ and it is locally compact, invariant and laminates into a collection of minimal sets. Moreover, the restriction of $\tau$ to $E$ is continuous and this set concentrates the pullback and forwards dynamics of the semiflow, in the sense that it contains  the pullback and forwards limits of all the trajectories in $\Omega \times \Int \wit K_B$, uniformly on bounded sets. Finally, when in addition the semiflow admits a point of strong sublinearity, we prove that it is dissipative in $\W\times \Int\wit K_B$,  $E$ is the global attractor and  it is a copy of the base, that is, $E=\{(\omega,b(\omega)) \mid \omega \in \Omega \}$.  This is the version in this context of the analogous result proved, when the interior of the positive cone is not empty, in N\'{u}\~{n}ez et al.~\cite{nuos2} and Zhao~\cite{zhao}.\smallskip

The aim of the second part of the paper is to  apply the theory just developed for delay systems which are monotone and sublinear for the exponential ordering  to analyse the long-term dynamics of
almost periodic compartmental Nicholson systems.    First of all, we analyse the relations of the coefficients which imply that the Nicholson systems induce a global monotone and sublinear semiflow for the  exponential ordering given by a diagonal matrix $B$. Under these conditions the semiflow is in fact strongly sublinear. Thus, assuming the uniform persistence of the systems, the existence of a global attractor in
$\W\times \Int \wit K_B$ which is a globally asymptotically stable copy of the base is proved. From here, using arguments of the theory of almost periodic functions  and a collection of linear changes of variables, we obtain  a new and more general reformulation of these conditions implying  the existence of a global attractor in $\W\times \Int X_+$, which is a globally asymptotically stable  copy of the base. These new conditions improve the non-autonomous versions of the inequalities given in~\cite{smth91}. The same kind of ideas can be used in the study of the existence of the so-called {\em special solutions\/} for non-autonomous FDEs. This theory originated in the 1960s in works by Ryabov~\cite{ryab} (see also Driver~\cite{driv}). In particular, we obtain new conditions that improve the non-autonomous version of the inequalities  given in Pituk~\cite{pitu} for the existence of special solutions for an almost periodic scalar Nicholson equation.\smallskip

The paper is organised in four sections. Section~\ref{sec-preli} contains the conditions of monotonicity for the exponential ordering, the definitions of semiequilibria, and their implications on the persistence and the dissipativity of the systems. In Section~\ref{sec-mon-sub}, the dynamical theory for monotone and sublinear semiflows in the positive cone for the exponential ordering is developed. Finally, the previous theory is applied in Section~\ref{sec-nicholson} to  almost periodic  compartmental Nicholson systems.
\section{Preliminaries on the exponential ordering and some results for monotone delay FDEs}\label{sec-preli}
Let $(\W,\sigma,\R)$ be a minimal flow over a compact metric space $\W$ and let $X=C([-r,0],\R^m)$ be the Banach space of continuous functions with the sup-norm $\|\cdot\|_\infty$, i.e., $\n{\varphi}_\infty=\sup_{s\in[-r,0]}\n{\varphi(s)}$ for $\varphi\in X$,  where $\n{\cdot}$ denotes the maximum norm on $\R^m$.   As usual, for a continuous map $y:[-r,\infty)\to \R^m$ and a time $t\geq 0$, $y_t$ denotes the map in $X$ defined by $y_t(s)=y(t+s)$, $s\in [-r,0]$.\par

Let us consider a family of finite-delay FDEs over $\W$,
\begin{equation}\label{fdefamily}
y'(t)=F(\w{\cdot}t,y_t)\,,\quad t\geq 0\,,\quad \text{for each}\;\w\in\W\,,
\end{equation}
defined by a
function $F\colon\Omega\times X \to\R^m$, $(\w,\varphi)\mapsto F(\w,\varphi)$ satisfying:\par\smallskip

\begin{enumerate}[label=\upshape(\textbf{F\arabic*}),series=delay_properties,leftmargin=27pt]
\item\label{F1} $F$ is continuous on $\W\times X$, of class $C^1$ with respect to the second variable, and for each  bounded set $X_1\subset X$, $F(\W\times X_1)$ is a bounded set in $\R^m$.
\end{enumerate}
\par\smallskip

Under this assumption, the standard theory of FDEs (see Hale and Verduyn Lunel~\cite{have}) assures that, for each $\w\in\W$ and $\varphi\in X$, the system~\eqref{fdefamily}$_\w$ admits  a locally defined  unique solution $y(t,\w,\varphi)$ with initial value $\varphi$, that is, $y(s,\w,\varphi)= \varphi(s)$ for each $s\in[-r,0]$. As a consequence, the family~\eqref{fdefamily} induces a local continuous skew-product semiflow
\begin{equation}\label{skewC}
  \begin{array}{cccl}
    \tau \colon&\mathcal{U}\subset \R^+\times\W\times X& \longrightarrow & \W\times X\\
         &\qquad (t,\w,\varphi) & \mapsto &(\w{\cdot}t, y_t(\w,\varphi))\,
  \end{array}
\end{equation}
which preserves the flow on the base $\W$ and has the so-called {\em semicocycle\/} property on the fiber. Namely, if $(\w,\varphi)\in \W\times X$, then for all $ t,s\geq 0$ for which the terms are  defined,
\begin{equation}\label{semicocycle}
y_{t+s}(\w,\varphi) = y_t(\w{\cdot}s,y_s(\w,\varphi))\,.
\end{equation}
\par
The theory of cooperative systems of delay differential equations is quite well developed (see, e.g., Smith~\cite{smit} and the references therein). Taking the standard cone of positive maps $X_+=\{\phi\in X\mid \phi(s)\geq 0 \; \text{for}\; s\in [-r,0]\}$ with nonempty interior  $\Int X_+=\{\phi\in X\mid \phi(s)\gg 0 \; \text{for}\; s\in [-r,0]\}$, $X$ is a strongly ordered Banach space. Note that we use the standard notation for $y\in \R^m$: $y\geq 0$ means that all components are nonnegative and $y\gg 0$ means that all components are positive.
The induced partial order relation on $X$ is then given by:
\begin{equation*}
\begin{split}
 \phi\le \psi \quad &\Longleftrightarrow \quad \psi-\phi\in X_+\,;\\
 \phi< \psi  \quad &\Longleftrightarrow \quad \psi-\phi\in X_+\;\text{ and }\;\phi\ne \psi\,;
\\  \phi\ll \psi \quad &\Longleftrightarrow \quad \psi-\phi\in \Int X_+\,.\qquad\quad\quad~
\end{split}
\end{equation*}
It is a well-known result that, if $F$ satisfies the so-called {\em quasimonotone condition}
\par\smallskip
\begin{enumerate}[label=\upshape(\textbf{Q}),leftmargin=27pt]
\item\label{Q} If $\phi\leq \psi$ and $\phi_i(0)=\psi_i(0)$ for some $i$, then $F_i(\w,\phi)\leq F_i(\w,\psi)$ for all $\w\in\W$\,,
\end{enumerate}
\par\smallskip\noindent
then $\tau$ is a monotone semiflow, that is, if $\phi\leq \psi$ and  $\w\in\W$, then $y_t(\w,\phi)\leq y_t(\w,\psi)$ for all $t\geq 0$ where both terms are defined.

We now introduce the exponential order in $X$, following~\cite{smth91}. We say that an $m\times m$ matrix $B=[b_{ij}]$ is {\it quasipositive\/} or  {\it cooperative\/} if all the off-diagonal entries are nonnegative, that is, $b_{ij}\geq 0$ whenever $i\not= j$. Given a cooperative matrix $B$ and  considering
the componentwise partial ordering on $\mathbb{R}^m$,  we introduce the normal positive
cone with empty interior in $X$,
\begin{equation*}
 K_B=\big\{ \phi\in X\mid \phi\geq 0 \;\text{ and }\; \phi(t)\geq
          e^{B(t-s)}\phi(s)\quad
          \text{for } -r\leq s\leq t\leq 0\big\}
\end{equation*}
which induces the following partial order relation on $X$, called the exponential ordering:
\begin{align*}
  \phi\le_B \psi\;&\Longleftrightarrow\; \phi\le \psi \;\text{ and }\;
              \psi(t)-\phi(t)\geq e^{B(t-s)}(\psi(s)-\phi(s))\,,\;
              -r\leq s\leq t\leq 0\,,\\
  \phi<_B \psi \;&\Longleftrightarrow\; \phi \leq_B \psi\;\;\text{ and }\;
             \phi\neq \psi\,.
\end{align*}
A smooth map $\phi$ belongs to $K_B$ if and only if $\phi\geq 0$ and $\phi'\geq B\,\phi$ on $[-r,0]$.\par\smallskip
\par
Let us assume one additional condition on $F$:\par\smallskip

\begin{enumerate}[resume*=delay_properties]
\item\label{F2}   $F(\w,\psi)-F(\w,\phi)\geq B\,(\psi(0)-\phi(0))\;$  whenever $\;\phi \leq_B\psi\;$  and  $\;\w\in\W$\,,
\end{enumerate}\par\smallskip
\noindent which is a necessary and sufficient condition for the monotonicity of the semiflow~\eqref{skewC} for the exponential ordering (see~\cite[Proposition~3.1]{smth91}), that is,  if $\phi\leq_B\psi$ and $\w\in\W$, then
\begin{equation*}
 y_t(\w,\phi)\leq_B y_t(\w,\psi) \quad \text{for all } t\geq 0\; \text{for which both are defined}.
\end{equation*}
Besides, as proved in~\cite[Corollary~3.2]{smth91}, a componentwise separating behaviour is also derived from~\ref{F2}.
Namely, if $\phi,\psi\in X$ satisfy $\phi\leq_B \psi$ and, for some $\w\in\W$, some component $1\leq i\leq m$, and some $t_0\geq -r$,  $y_i(t_0,\w,\phi)<y_i(t_0,\w,\psi)$, then $y_i(t,\w,\phi)<y_i(t,\w,\psi)$ for all $t\geq t_0$ for which both are defined.\par\smallskip
A stronger monotonicity condition presented in~\cite{smth91} is the following:\par\smallskip
\begin{enumerate}[resume*=delay_properties]
\item\label{F3}
$F(\w,\psi)-F(\w,\phi)\gg B\, (\psi(0)-\phi(0))\,$  whenever $\,\phi \leq_B\,\psi$, $\phi\ll\psi$   and $\w\in\W\,.$
\end{enumerate}
\par\smallskip
Condition~\ref{F3} implies \ref{F2}  by a continuity argument.  With this condition, the monotonicity of $\tau$ becomes slightly stronger, in a very precise sense. As stated above, the interior of the cone $K_B$ is empty in $X$. However, by taking the Banach space of Lipschitz functions $X_L\subset X$, with the usual Lipschitz norm
\[
\n{\phi}_L=\n{\phi}_\infty+\sup \left\{ \left|\frac{\phi(t)-\phi(s)}{t-s}\right|, \; s\neq t,\;s,t \in[-r,0] \right\},
\]
the restriction of the cone  $K_B$ to $X_L$ provides us with the cone of positive maps
\begin{equation*}
\wit K_B=\{\phi\in X_L\mid \phi\geq 0 \;\text{\;and}\; \phi'-B\,\phi\geq 0 \,\;\text{a.e. on}\; [-r, 0]\}
\end{equation*}
which instead has a nonempty interior in $X_L$:
\begin{align*}
  \Int\wit K_B & =\{\phi\in X_L\mid \phi\gg 0 \;\;\text{and}\; \phi'-B\,\phi\geq \bar\varepsilon \,\;\text{a.e. on}\; [-r, 0],\;\text{for some}\;\varepsilon>0\} \\
   & =\{\phi\in X_L\mid  \phi(-r)\gg 0 \text{ and  ess inf}\,(\phi'- B\,\phi)_i>0,\; 1\leq i\leq m\}\,.
\end{align*}
Here, and all through the paper, $\bar \varepsilon$ stands for the vector (or the constant map) whose components identically equal the value $\varepsilon$.
Now we can write $\phi \gg_B 0$ meaning that $\phi \in \Int\wit K_B$ and
\[
\phi\ll_B \psi \;\Longleftrightarrow\; \psi - \phi \in \Int\wit K_B \,.
\]
Note that the Lipschitz character of the maps involved is implicit whenever the strong order $\ll_B$ is used.

At this point one may wonder whether we could not just work on the space $X_L$, but, even if the section $\tau_t:\W\times X_L \to \W\times X_L$ is continuous for each fixed $t\geq 0$,  the semiflow $\tau:\R^+\times \W\times X_L \to \W\times X_L$ is not continuous for $t$ in the range of times  $[0,r]$: see~\cite[pp.~345]{smth91} for further explanations. What is known is that
\begin{equation}\label{XXL}
  \begin{array}{cccl}
    \tau\colon &[r,\infty)\times \W\times X& \longrightarrow & \W\times X_L\\
         &\;(t,\w,\varphi) & \mapsto &(\w{\cdot}t, y_t(\w,\varphi))
  \end{array}
\end{equation}
is continuous, as far as the solution $y_t(\w,\varphi)$ is defined. For these reasons, both spaces are going to interplay.

These are the strong monotonicity relations we get when condition~\ref{F3} holds.
\begin{prop}\label{prop monot fuerte}
Assume that conditions~\ref{F1} and~\ref{F3} hold. Then:
\begin{itemize}
 \item[\rm{(i)}]  If  $\phi\leq_B\psi$, $\phi(0)\ll\psi(0)$ and $\w\in\W$, then
 $y_t(\w,\phi)\ll_B y_t(\w,\psi)$  for each $t\geq 2r$ for which both  are defined.
 \item[\rm{(ii)}]  If
$\phi\ll_B\psi$, then $y_t(\w,\phi)\ll_B y_t(\w,\psi)$  for each $t\geq 0$ for which both  are defined.
\end{itemize}
\end{prop}
\begin{proof}
Recall that condition~\ref{F3} implies condition~\ref{F2}. As mentioned above, thanks to condition~\ref{F2}, if $\phi\leq_B\psi$ with $\phi(0)\ll\psi(0)$, then $y(t,\w,\phi)\ll y(t,\w,\psi)$ for all $t\geq 0$ for which both are defined, so that $y_t(\w,\phi)\ll y_t(\w,\psi)$ for all $t\geq r$.  Note that $y(t,\w,\phi)$ and $y(t,\w,\psi)$ are both differentiable, and in particular locally Lipschitz for $t>0$. Then, by condition~\ref{F3}, for all $t\geq r$,
\begin{multline}\label{eq:aux}
 y'(t,\w,\psi) -y'(t,\w,\phi)-B\,(y(t,\w,\psi) -y(t,\w,\phi)) \\
 = F(\w{\cdot}t,y_t(\w,\psi))-F(\w{\cdot}t,y_t(\w,\phi))-B\,(y(t,\w,\psi) -y(t,\w,\phi))\gg 0\,.
\end{multline}
By continuity, for each $t\geq 2r$ we can take an $\varepsilon>0$ (which depends on $t$) such that  $y'(t+s,\w,\psi) -y'(t+s,\w,\phi)-B\,(y(t+s,\w,\psi) -y(t+s,\w,\phi))\geq \bar \varepsilon$ for all $s\in [-r,0]$. All in all, $y_t(\w,\psi)-y_t(\w,\phi)\gg_B 0$ for all $t\geq 2r$ for which both are defined, as asserted in~(i).

As for~(ii), if $\phi\ll_B\psi$, we already know that $y_t(\w,\phi)\leq_B y_t(\w,\psi)$ and, since in particular $\phi\ll \psi$, also  $y_t(\w,\phi)\ll y_t(\w,\psi)$  for each $t\geq 0$ for which both solutions are defined. Then, relation~\eqref{eq:aux} holds for all $t\geq 0$. From this and the fact that $\phi\ll_B\psi$, the result follows easily. The proof is finished.
\end{proof}
In monotone non-autonomous dynamical systems, the so-called semiequilibria, introduced both for the deterministic and random cases (see~\cite{nono2} and~\cite{book:chue}),  are useful objects to determine invariant zones (see also~Novo and Obaya~\cite{noob1}, among others). We recall the definitions for the exponential ordering, which were introduced by Novo et al.~\cite{noov} in a context of infinite-delay neutral FDEs.
\begin{definition} A   map  $a\colon \W \to X$ such that $y_t(\w,a(\w))$ exists for all $t\geq 0$  is
\begin{itemize}
\item a \emph{sub-equilibrium for the exponential ordering}  if $a(\w{\cdot}t)\leq_B y_t(\w,a(\w))$ for each $\w\in\W$ and $t\geq 0\,$;
\item a \emph{super-equilibrium for the exponential ordering}  if $a(\w{\cdot}t)\geq_B y_t(\w,a(\w))$ for each $\w\in\W$ and $t\geq 0\,$;
\item an \emph{equilibrium} if $a(\w{\cdot}t)= y_t(\w,a(\w))$ for each $\w\in\W$ and $t\geq 0\,$.
\end{itemize}
 A sub-equilibrium (resp. super-equilibrium) for the exponential ordering  is said to be \emph{strong} if there exists a time $s_\ast>0$ such that $a(\w{\cdot}s_\ast)\ll_By_{s_\ast}(\w, a(\w))$ (resp. $a(\w{\cdot}s_\ast)\gg_By_{s_\ast}(\w, a(\w))$) for each $\w\in\W$.
\end{definition}
The next result shows how sub-equilibria and strong sub-equilibria  for $\leq_B$ can be constructed from a family of sub-solutions.  By changing the sign of the inequalities,  we obtain the corresponding results for super-equilibria and strong super-equilibria in terms of a family of super-solutions. The following definition is needed.
\begin{definition}
A map  $\wit a\colon \W\to\R^m$  is said to be \emph{$C^1$\,along the trajectories of the base flow}  if the map $\wit a_\w\colon\R\to\R^m$, $t\mapsto \wit a(\w{\cdot}t)$ is of class $C^1$ for each $\w\in\W$. We will denote $\wit a'(\w)= (\wit a_\w)'(0)=(d/dt) \wit a(\w{\cdot}t)|_{t=0}$.
\end{definition}
We remark that in many occasions one deals with semicontinuous semiequilibria. This permits the combination of both dynamical and topological properties, with important consequences. The reader is referred to~\cite{nono2} for a precise statement of these properties. In particular, a semicontinuous semiequilibrium has a residual set of points of continuity.
\begin{prop}\label{prop:subequi}
Assume  conditions~\ref{F1} and \ref{F2} and let $\wit a\colon \W\to\R^m$ be $C^1$\,along the trajectories of the base flow. Consider the map $a\colon\W\to X$, $\w\mapsto a(\w)$ defined by  $a(\w)(s)=\wit a(\w{\cdot}s)$ for each $s\in[-r,0]$. Then:
\begin{itemize}
\item[\rm{(i)}] The map $a$
  is a sub-equilibrium for the exponential ordering provided that for each $\w\in\W$, $y(t,\w,a(\w))$ is defined for all  $t\geq 0$ and
\begin{equation}\label{subsolution}
\wit a'(\w)\leq F(\w,a(\w))\,.
\end{equation}
\item[\rm{(ii)}]  If the conditions in \rm{(i)} and \ref{F3} hold, and there is an $\w_0\in\W$ which is simultaneously a continuity point of $a$ and $a\circ \sigma_{s_\ast}$, for some $s_\ast\geq 3r$,  such that
\begin{equation}\label{stronsubequi}
\wit a'(\w_0)\ll F(\w_0,a(\w_0))\,,
\end{equation}
then $a$ is a strong sub-equilibrium for the exponential ordering.
\end{itemize}
\end{prop}
\begin{proof} (i) The proof is inspired in that of~\cite[Proposition~3.1]{smth91}. Let $F^\ep=F + \bar \ep\,$ for $\ep>0$ and let $y^\ep(t,\w,a(\w))$ denote the solution of $y'(t)=F^\ep(\w{\cdot}t,y_t)$ with initial value $a(\w)$. We claim that
$y_t^\ep(\w,a(\w))\geq_B  a(\w{\cdot}t)$ for each $\w\in\W$, $t\geq 0$ and $\ep>0$, and hence the proof is finished by letting $\ep \downarrow 0$.
\par
In order to see it, we consider $x^\ep(t)=y^\ep(t,\w,a(\w))- \wit a_\w(t)$, which satisfies $x^\ep(s)=0$ for each $s\in[-r,0]$, and  the closed set $I=\{t\in[0,\infty)\mid x_t^\ep\geq_B 0\}$. Notice that $0\in I$.  Note also that we are implicitly assuming that $y^\ep(t,\w,a(\w))$ is extendable to $[0,\infty)$. If not, we work on its maximal interval of existence, but recall that, given any compact set $[0,T]$, $y^\ep(t,\w,a(\w))$ is defined on $[0,T]$ for every $\ep$ small enough, since $y(t,\w,a(\w))$ is assumed to be globally defined and $F^\ep=F + \bar \ep\,$.
 Next we check that given $t\in I$ there is a $\delta>0$ such that $[t,t+\delta)\subset I$.
Inequality~\eqref{subsolution} can be written as $(\wit a_\w)'(t)\leq F(\w{\cdot}t,(\wit a_\w)_t)$ for each $t\geq 0$ and $\w\in \W$. Then, from condition~\ref{F2} and $t\in I$, we deduce that
\begin{align*}
 (x^\ep)'(t)- B\, x^\ep(t)& =F(\w{\cdot}t,y^\ep_t(\w, a(\w)))+\bar\ep- (\wit a_\w)'(t)- B\, x^\ep(t)\\
    & \geq \bar\ep + F(\w{\cdot}t,y^\ep_t(\w, a(\w)))- F(\w{\cdot}t,a(\w{\cdot}t)) - B\, x^\ep(t)\geq \bar\ep\gg 0\,.
\end{align*}
Hence, there is a $\delta>0$ such that   $(x^\ep)'- B\, x^\ep\geq 0$  on $[t,t+\delta)$. Note that in particular $x^\ep(t)\geq 0$ and $(x^\ep)'\geq B\, x^\ep$  on $[t,t+\delta)$. Since $B$ is a cooperative matrix, we can deduce by comparing the solutions that $x^\ep(s)\geq 0$ for all $s\in[t,t+\delta)$. As a result, $x_s^\ep\geq_B 0$ for all $s\in[t,t+\delta)$, that is, $[t,t+\delta) \subset I$, as wanted. Therefore, $I=[0,\infty)$, which finishes the proof of~(i).
\par
(ii) From~(i) we know that $a$ is a sub-equilibrium for the exponential ordering, that is, $a(\w{\cdot}t)\leq_B y_t(\w,a(\w))$ for each $\w\in\W$ and  $t\geq 0$. Therefore, according to~\cite[Proposition~4.2(i)]{nono2}, which can be easily generalised to this case, to show that $a$ is a strong sub-equilibrium for the exponential ordering, it is enough to check the existence of  a time $s_\ast>0$  and a point $\w_0\in\W$ which is simultaneously a continuity point of $a$ and $a\circ \sigma_{s_\ast}$ such that $a(\w_0{\cdot}s_\ast)\ll_B y_{s_\ast}(\w_0,a(\w_0))$. \par In order to check this, notice that condition~\eqref{stronsubequi}  and the continuity of the maps involved imply the existence of an $\ep\in (0,r)$ such that $(\wit a_{\w_0})'(t)\ll F(\w{\cdot}t,(\wit a_{\w_0})_t)$ for each $t\in[0,\ep)$, and we can take a time $t_0\in (0,\ep)$ satisfying $\wit a_{\w_0}(t_0)\ll y(t_0,\w_0,a(\w_0))$, i.e., $a(\w_0{\cdot} t_0)(0)\ll y_{t_0}(\w_0,a(\w_0))(0)$. Thus, since $a(\w_0{\cdot} t_0)\leq_B y_{t_0}(\w_0,a(\w_0))$, by Proposition~\ref{prop monot fuerte}(i),
\[
y_t(\w_0{\cdot}t_0, a(\w_0{\cdot}t_0)) \ll_By_t(\w_0{\cdot}t_0,y_{t_0}(\w_0,a(\w_0))) \quad \text{for each } t\geq 2r\,,
\]
and then,  $a(\w_0{\cdot}(t+t_0))\ll _B y_{t+t_0}(\w_0,a(\w_0))$. Thus,  choosing $t=s_\ast-t_0\geq 2r$, we conclude that  $a(\w_0{\cdot}s_\ast)\ll_B  y_{s_\ast}(\w_0,a(\w_0))$,  which finishes the proof.
\end{proof}

As a consequence, we determine conditions concerning the permanence of the solutions in the positive cone $K_B$ and in the interior of the positive cone $\Int\wit K_B$.
\begin{coro}\label{coro:invariance cone}
Assume \ref{F1}--\ref{F2} and $F(\w,0)\geq 0$ for each $\w\in\W$.
Then:
\begin{itemize}
\item[\rm{(i)}]  $y_t(\w,\psi)\geq _B 0$ for each $\w\in\W$, $\psi\geq_B 0$ and $t\geq 0$ in the domain of definition.
\item[\rm{(ii)}] If \ref{F3} also holds, then $y_t(\w,\psi)\gg _B 0$ for each  $\w\in\W$, $\psi\gg_B 0$ and $t\geq 0$ in the domain of definition.
\end{itemize}
\end{coro}
\begin{proof}
Note that Proposition~\ref{prop:subequi}(i) implies that $y_t(\w,0)\geq_B 0$ as long as $y(t,\w,0)$ is defined.   Then, by monotonicity,  we deduce that, if   $\psi\geq _B0$, then $y_t(\w,\psi)\geq_B y_t(\w,0)\geq_B 0$ for all $t\geq 0$ in the domain of definition of $\psi$, so that~(i) holds.

For~(ii), take $\w\in\W$ and $\psi\gg_B 0$. Then, by Proposition~\ref{prop monot fuerte}(ii), $y_t(\w,\psi)\gg _B y_t(\w,0)\geq_B 0$, whenever defined.  The proof is finished.
\end{proof}

In the rest of this section, we take advantage of    the presence of semiequilibria. When no mention is made,  we will be assuming that the semiflow is globally defined. First, we consider the dynamical property of uniform persistence. We include the standard  definitions for the usual ordering and for the exponential ordering, respectively. We remark that these two properties are in principle unrelated.
\begin{definition}\label{defi persistence hull}
(i) The skew-product semiflow $\tau$ induced by the  family of systems~\eqref{fdefamily} is {\em uniformly persistent} in the interior of the positive cone $\Int X_+$ if there is a map $\psi\gg 0$ such that, for every $\w\in \W$ and every initial map $\varphi\gg 0$,  there exists a time $t_0=t_0(\w,\varphi)$ such that      $y_t(\w,\varphi)\geq \psi$ for all $t\geq t_0$.

(ii) The skew-product semiflow $\tau$ induced by the  family of systems~\eqref{fdefamily} is {\em uniformly persistent} in the interior of the positive cone $\Int \wit K_B$ associated to the exponential ordering $\leq_B$ if there is a map $\psi\gg_B 0$ such that, for every $\w\in \W$ and every initial map $\varphi\gg_B 0$,  there exists a time $t_0=t_0(\w,\varphi)$ such that
     $y_t(\w,\varphi)\geq_B \psi$ for all $t\geq t_0$.
\end{definition}
We give an easy-to-check sufficient condition for the uniform persistence in the interior of the positive cone $\Int \wit K_B$.
\begin{prop}\label{prop:unif pers exp ord}
Assume conditions~\ref{F1} and~\ref{F3} hold, and $F(\w,0)\geq 0$ for each $\w\in\W$. If there exists an $\w_0\in\W$ such that $F(\w_0,0)\gg 0$, $\tau$ is uniformly persistent in the interior of the positive cone $\Int \wit K_B$.
\end{prop}
\begin{proof}
From Proposition~\ref{prop:subequi}(ii) we deduce that the null map from $\W$ to $X$ is a continuous strong sub-equilibrium for $\leq_B$ and, reasoning as in~\cite[Proposition~4.3]{nono2} for the exponential ordering, there exist a $\psi\gg_B 0$ and a time $s>0$ such that $a_s(\w):=y_{s}(\w{\cdot}(-s),0)\geq_B \psi$ for each $\w\in\W$. In particular, since for all $\w$, $a_t(\w)\geq_B a_s(\w)$ whenever  $t\geq s$, evaluating at $\w{\cdot}t$ we obtain that $y_t(\w,0)\geq_B \psi$ for each $\w\in\W$ and $t\geq s$. Finally, from Proposition~\ref{prop monot fuerte}(ii), if $\varphi\gg_B 0$, we deduce that $y_t(\w,\varphi)\gg_B y_t(\w,0)$ for each $t\geq 0$. Hence, $y_t(\w,\varphi)\geq_B \psi$ for each $\w\in\W$ and $t\geq s$, which shows the claimed uniform persistence.
\end{proof}
Finally, we concentrate on the existence of attractors. As in the case of standard dynamical systems, also in non-autonomous dynamical systems the existence of a global attractor permits a better understanding of the dynamics (see, e.g.,~\cite{klra}). A nice reference for processes and pullback attractors, which are crucial in non-autonomous dynamics, is~\cite{book:CLR}. For a skew-product semiflow over a compact base $\W$, {\em the global attractor\/}   $\A\subset \W\times X$, when it exists,  is an invariant compact set attracting a certain class of sets in $\W\times X$ forwards in time; namely,
\[
\lim_{t\to\infty} {\rm dist}(\tau_t(\W\times X_1),\A)=0 \quad\text{for each}\; X_1\in \mathcal{D}(X)\,,
\]
for the Hausdorff semidistance ${\rm dist}$. The standard  choices for $\mathcal{D}(X)$ are either the class of bounded subsets $\mathcal{D}_b(X)$ of $X$ or that of compact subsets $\mathcal{D}_c(X)$ (see~\cite{chks}).

For the sake of completeness, we recall the definition of the Hausdorff semidistance between subsets of a metric space $(Y,d)$. Given $Y_1,Y_2\subset Y$,
\begin{equation}\label{semidist}
  {\rm dist}(Y_1,Y_2) := \sup_{y_1\in Y_1}\inf_{y_2\in Y_2} d(y_1,y_2)=\sup_{y_1\in Y_1} {\rm d}(y_1,Y_2)\,,
\end{equation}
where ${\rm d}(y_1,Y_2)$ denotes the usual distance between the point $y_1$ and the set $Y_2$. The Hausdorff metric between two compact sets of $Y$ is then defined by
\[
{\rm dist_{\mathcal H}}(Y_1,Y_2) := \max\{{\rm dist}(Y_1,Y_2),{\rm dist}(Y_2,Y_1)\}\,.
\]

As $\W$ is compact,  the non-autonomous set $\{A(\w)\}_{\w\in \W}$, formed by the $\w$-sections of $\A$ defined by $A(\w)=\{\varphi\in X\mid (\w,\varphi)\in \A\}$ for each $\w\in \W$, is a {\em pullback attractor\/}, that is,  $\{A(\w)\}_{\w\in \W}$ is compact, invariant, and it pullback attracts all the sets $X_1 \in \mathcal{D}(X)$:
\begin{equation*}
\lim_{t\to\infty} {\rm dist}(y_t(\w{\cdot}(-t),X_1),A(\w))=0\quad \text{for all}\; \w\in \W\,.
\end{equation*}
The implications for each fixed $\w\in \W$ are expressed in terms of {\em pullback attraction for the related evolution process on $X$\/}, defined by $S_{\w}(t,s)\,\varphi=y_{t-s}(\w{\cdot}s,\varphi)$ for each $\varphi\in X$ and $t\geq s$. Precisely,
for each fixed $\w\in \W$, the family of compact sets $\{A(\w{\cdot}t)\}_{t\in  \R}$ is the pullback attractor for the process  $S_{\w}(\cdot,\cdot)$,
 meaning that:
\begin{itemize}
\item[(i)] it is invariant, i.e.,  $S_{\w}(t,s)\,A(\w{\cdot}s)=A(\w{\cdot}t)$ for all $t\geq s\,$;
\item[(ii)] it pullback attracts the class of sets $\mathcal{D}(X)$, i.e., for each  $X_1\in \mathcal{D}(X)$,
\[
\lim_{s\to -\infty}{\rm dist}(S_{\w}(t,s)\,X_1,A(\w{\cdot}t))=0 \quad \hbox{for all}\; t\in \R\,;
\]
\item[(iii)] it is the minimal family of closed subsets of $X$ with property~(ii).
\end{itemize}
However, the forwards and pullback dynamics are in general unrelated, so that it is always an interesting  problem to know whether the pullback attractor is also a {\em forwards attractor\/} for the process. That is to know if, for each $\w\in\W$,
\begin{equation*}
\lim_{t\to \infty}{\rm dist}(y_t(\w,X_1),A(\w{\cdot}t))=0   \quad \hbox{for each}\; X_1\in \mathcal{D}(X)\,.
\end{equation*}

We finish this section with a result on the existence of a global attractor in the positive cone and in the interior of the positive cone, in terms of the existence of continuous strong sub and super-equilibria, for the exponential ordering. Note that this covers a huge range of applications to delay systems in biology or ecology, where only positive solutions make sense. We remark that we impose a restriction on the class of quasipositive matrices $B$ defining the ordering $\leq_B$, which is satisfied by the commonly used diagonal matrices with negative diagonal elements.
\begin{theorem}\label{teor:atractores}
Assume that the quasipositive matrix $B=[b_{ij}]$ defining the exponential ordering satisfies that $b_{ii}< -\sum_{j\not=i}b_{ij}$,   conditions~\ref{F1} and~\ref{F3} hold, and $F(\w,0)\geq 0$ for each $\w\in\W$. Then:
\begin{itemize}
\item[{\rm (i)}]   If there exist an $R_0>0$ and an $\w_0\in\W$ such that $F(\w,\bar R)\leq 0$ and $F(\w_0,\bar R)\ll 0$ for each  $\w\in \W$ and $R\geq R_0$, then the semiflow $\tau$ is globally defined on $\W\times K_B$ and there exists a global attractor  with respect to the class $\mathcal{D}_b(K_B)$ of bounded subsets in the positive cone.
\item[{\rm (ii)}] If the condition in {\rm(i)} holds and, besides, there exists an $\w_1\in \W$ such that $F(\w_1,0)\gg 0$, then there exists a global attractor in $\W\times \Int \wit K_B$ with respect to the class $\mathcal{D}_c(\Int \wit K_B)$ of compact sets in $X_L$ contained in $\Int \wit K_B$.
\end{itemize}
\end{theorem}
\begin{proof}
First of all, Corollary~\ref{coro:invariance cone} gives us the invariance of both the cone of positive elements and its interior, that is to say, of the sets $\W\times K_B$ and $\W\times \Int \wit K_B$, respectively. Also note that the condition imposed on $B$ guarantees that the constant map $\bar R\in \Int \wit K_B$ for all $R> 0$.\smallskip

(i) Notice that Proposition~\ref{prop:subequi}(ii) implies that the constant maps $\W\to K_B$, $\w\mapsto \bar R$ are continuous strong super-equilibria for the exponential ordering for each $R\geq R_0$. In particular this means that the sets $\{(\w,\phi)\in \W\times K_B\mid \phi\leq_B \bar R\}$ are positively invariant. Since, thanks to the restriction imposed on $B$, given a $\phi\in K_B$ we can take a sufficiently big $R>0$ so that $\phi\leq_B \bar R$, we can deduce that  all the solutions are bounded and the semiflow is globally defined on $\W\times K_B$.\smallskip

To get the existence of a global attractor, it suffices to find an  absorbing compact set (see, e.g., \cite[Theorem~1.36]{klra}). In our skew-product setting, we search for a compact set $H\subset K_B$ such that for each bounded set $X_1\subset K_B$ there exists a $t_1=t_1(X_1)$ such that $\tau_t(\W\times X_1)\subset  \W\times H$ for all $t\geq t_1$. Thanks to condition~\ref{F1}, the map $y_r:\W\times X\to X$ is compact, meaning that it takes bounded sets into relatively compact sets. Then, it is immediate that the set
\[
H=\cls \big\{y_r(\w,\varphi) \mid \w\in \W\,,\; 0 \leq_B \varphi \leq_B \bar R_0 \big\}
\]
is compact in $K_B$. Let us check that the compact set $\W\times H$ is absorbing. For each bounded set $X_1\subset K_B$, we can find an $R\geq R_0$ big enough so that $0\leq_B \phi \leq_B \bar R$ for all $\phi\in X_1$ and, by monotonicity, $0\leq_B y_t(\w,\phi) \leq_B y_t(\w,\bar R)$ for all  $\w\in\W$, $\phi\in X_1$, and $t\geq 0$.  Since the constant maps $\bar R$ are strong super-equilibria for $\leq_B$ for all $R\geq R_0$, it is not difficult to deduce that there exists a $t_0\geq 0$ such that
$0\leq_B y_t(\w,\phi) \leq_B \bar R_0$ for all $(\w,\phi)\in\W\times X_1$ and $t\geq t_0$. Now, it suffices to apply the semicocycle property~\eqref{semicocycle} to get that, if $t\geq t_1:= t_0+r$, then $y_t(\w,\phi)=y_r(\w{\cdot}(t-r),y_{t-r}(\w,\phi))\in H$ for all $(\w,\phi)\in\W\times X_1$, as we wanted to see.\smallskip

(ii) First of all, recall that the interior of the positive cone $K_B$ is empty in $X$ but it is nonempty when restricted to $X_L$.  For this reason we consider compact sets in $X_L$ when we look for a global attractor in $\W\times\Int\wit K_B$.\smallskip

Since $F(\w_1,0)\gg 0$, $\tau$ is uniformly persistent in $\Int\wit K_B$ by  Proposition~\ref{prop:unif pers exp ord}.  Let  $\psi\gg_B 0$ be the map in Definition~\ref{defi persistence hull}(ii) and  consider the compact set
$H=\cls \big\{y_r(\w,\varphi) \mid \w\in \W\,,\; \psi \leq_B \varphi \leq_B \bar R_0 \big\}$ in $K_B$.
Then,  take $\wit H=\{y_r(\w,\varphi) \mid (\w,\varphi)\in \W\times H\}$. $\wit H$ is a compact set in $X_L$ as it is  the image by the continuous map $y_r$ of a compact set in $\W\times X$ (see~\eqref{XXL}). It is easy to check that $\wit H\subset \Int \wit K_B$ because the map $\w\in\W\mapsto y_r(\w,\psi_0)\in \Int\wit K_B\subset X_L$ is continuous for each $\psi_0\gg_B 0$ and $\W$ is compact. Let us see that $\W\times \wit H$ is absorbing. Given a compact set $X_1\subset \Int \wit K_B$, there exist a $\varphi_0\gg_B 0$ and an $R\geq R_0$ big enough such that  $\varphi_0\leq_B \phi \leq_B  \bar R$ for every $\phi\in X_1$. By monotonicity, $y_t(\w,\varphi_0)\leq_B y_t(\w,\phi) \leq_B y_t(\w,\bar R)$ for all $(\w,\phi)\in\W\times X_1$ and $t\geq 0$. In fact, as seen in the proof of  Proposition~\ref{prop:unif pers exp ord}, for every $\varphi_0\gg_B 0$ we get a time $t_0=t_0(\varphi_0)$ such that $\psi\leq_B y_t(\w,\varphi_0)$ for all $t\geq t_0$ and $\w\in\W$, As in (i), by taking $t_0$ bigger if necessary, we can assert that  $\psi\leq_B y_t(\w,\phi)\leq_B \bar R_0$ for all $(\w,\phi)\in\W\times X_1$ and $t\geq t_0$. Now, by writing
$y_t(\w,\phi)=y_r(\w{\cdot}(t-r),y_r(\w{\cdot}(t-2r),y_{t-2r}(\w,\phi)))$ it is easy to see that $y_t(\w,\phi)\in \wit H$ for all $(\w,\phi)\in\W\times X_1$ and $t\geq t_1:=t_0+2r$. The proof is finished.
\end{proof}
\section{Long-term behaviour in the monotone and sublinear case}\label{sec-mon-sub}
In this section we consider a family of FDEs~\eqref{fdefamily} satisfying \ref{F1} and \ref{F2}, so that the  induced skew-product semiflow is monotone for the exponential ordering.
As explained in the previous section,  the skew-product semiflow $\tau:\R^+\times \W\times X_L  \to \W\times X_L $ is not continuous for $t\in [0,r]$, whereas it is continuous for $t\in [r,\infty)$.
\par
By adding some extra conditions on $F$, including a sublinear condition which results in the invariance of the positive cone $K_B$ and the sublinerity of the semiflow for the exponential ordering, we are able to determine a special subset in $\W\times \Int \wit K_B\subset \W\times X_L$ over which the restriction of the semiflow  is continuous and it concentrates the essential information for both the long-term behaviour of the solutions in $\Int \wit K_B$, and for attraction.
\par\smallskip
Namely, in this section we assume the following sublinearity condition:
\par\smallskip
\begin{enumerate}[resume*=delay_properties]
\item\label{F4}
$F(\w,\lambda\, \psi)\geq \lambda\, F(\w,\psi)\,$  whenever $\,\psi\geq _B 0$\,,\; $\lambda \in [0,1]\,,$\;   and $\,\w\in\W\,.$
\end{enumerate}
\par\smallskip\noindent
Next we prove how the properties on the systems are transferred to their solutions.
\begin{prop}\label{prop:sublinear}
Assume that conditions \ref{F1}, \ref{F2}, and \ref{F4} hold. Then, the semiflow is globally defined on $\W\times K_B$ and it is sublinear for the exponential ordering, that is,  for each  $\w\in\W$, $\psi\geq_B 0$, and $\lambda \in[0,1]$,
\begin{equation*}
 y_t(\w,\lambda\,\psi)\geq_B \lambda\, y_t(\w,\psi) \quad\text{for all } t\geq 0\,.
\end{equation*}
\end{prop}
\begin{proof}
Thanks to condition \ref{F4}, in particular $F(\w,0)\geq 0$ for all $\w\in\W$ so that, if $\w\in\W$ and  $\psi\geq_B0$, then  $y_t(\w,\psi)\geq_B 0$ whenever defined, by Corollary~\ref{coro:invariance cone}(i). Once we know that the solutions starting in $\W\times K_B$ remain in $\W\times K_B$ while defined, we only worry about boundedness above. Let us check that there exists a $c>0$ such that $F(\w,\psi)\leq (1+\|\psi\|_\infty)\,\bar c$ for all $(\w,\psi)\in \W\times K_B$, for the  vector $\bar c\in\R^m$. To see it,  by~\ref{F1}, if $\|\psi\|_\infty\leq 1$, then we can choose a  $c>0$ big enough so that  $F(\w,\psi)\leq \bar c$, for all $\w\in\W$. And, if $\|\psi\|_\infty>1 $, then we can deduce from~\ref{F4} that
\[
F(\w,\psi)=F\left(\w,  \|\psi\|_\infty\,\frac{\psi}{\|\psi\|_\infty}\right)\leq \|\psi\|_\infty\,F\left(\w,  \frac{\psi}{\|\psi\|_\infty}\right)\leq \|\psi\|_\infty\,\bar c\,,
\]
so that we are done.
Now note that the solutions of the family of systems $y'(t)=\wit F(\w{\cdot}t,y_t)$, $\w\in\W$, given by the globally Lipschitz map in the second variable $\wit F(\w,\psi)=(1+\|\psi\|_\infty)\,\bar c$, are globally defined, and that $\wit F$ satisfies the quasimonotone condition \ref{Q} when restricted to $\W\times K_B$. Therefore, we can compare solutions  (see~\cite[Theorem~5.1.1]{smit}) to deduce that $y(t,\w,\psi)$ is also globally defined for all  $(\w,\psi)\in \W\times K_B$.

Next we check the sublinearty for fixed $\w\in\W$ and $\psi\geq_B 0$. Since $y_t(\w,\psi)\geq_B 0$ for all $t\geq 0$, from~\ref{F4} we have that
\begin{equation}\label{eq F}
F(\w{\cdot}t,\lambda\, y_t(\w,\psi))\geq \lambda\, F(\w{\cdot}t,y_t(\w,\psi))\quad\text{for each } \lambda\in[0,1] \text{ and }t\geq 0\,.
\end{equation}
The proof now continues as that of Proposition~\ref{prop:subequi}.  We consider the function $F^\ep=F + \bar \ep\,$ for $\ep>0$ and let $y^\ep(t,\w,\varphi)$ denote the solution of $y'(t)=F^\ep(\w{\cdot}t,y_t)$ with initial value $\varphi$. We claim that
$y_t^\ep(\w,\lambda\,\psi)\geq_B \lambda\, y_t(\w,\psi)$ for all $t\geq 0$ and then the proof is finished by letting $\ep \downarrow 0$.
\par
In order to prove the claim, we take $x^\ep(t)=y^\ep(t,\w,\lambda\, \psi)-\lambda\, y(t,\w,\psi)$, which satisfies $x^\ep(s)=\lambda\,\psi(s)-\lambda\,\psi(s)=0$ for each $s\in[-r,0]$, and we consider the closed set $I=\{t\in[0,\infty)\mid x_t^\ep\geq_B 0\}$. Notice that $0\in I$. Next we check that given $t\in I$ there is a $\delta>0$ such that $[t,t+\delta)\subset I$. Since
\[
 (x^\ep)'(t)-B\, x^\ep(t)=F(\w{\cdot}t,y^\ep_t(\w,\lambda\,\psi))+\bar\ep-\lambda\, F(\w{\cdot}t,y_t(\w,\psi))-B\, x^\ep(t)\,,
\]
we can apply \eqref{eq F} to get
\begin{align*}
  (x^\ep)'(t)-B\, x^\ep(t)&\geq F(\w{\cdot}t,y^\ep_t(\w,\lambda\,\psi))-F(\w{\cdot}t,\lambda\, y_t(\w,\psi))\\
    & \qquad\quad -B\,(y^\ep(t,\w,\lambda\, \psi)-\lambda\, y(t,\w,\psi))+\bar\ep\,.
\end{align*}
Now, since $t\in I$ we know that $y^\ep_t(\w,\lambda\,\psi)\geq_B \lambda\, y_t(\w,\psi)$, and from~\ref{F2} we conclude that $(x^\ep)'(t)-B\, x^\ep(t)\geq \bar\ep$. Hence, there is a $\delta>0$ such that   $(x^\ep)'(s)-B\, x^\ep(s)\geq 0$  for $s\in [t,t+\delta)$. Besides, since in particular $x^\ep(t)\geq 0$ and $(x^\ep)'\geq B\, x^\ep$  on $[t,t+\delta)$ and $B$ is a cooperative matrix, by comparing the solutions we get that $x^\ep(s)\geq 0$ for $s\in[t,t+\delta)$. Consequently, $x_s^\ep\geq_B 0$ for each $s\in[t,t+\delta)$ and so $[t,t+\delta) \subset I$. Therefore, $I=[0,\infty)$ and the proof is finished.
\end{proof}
Hereafter, the existence of a bounded solution in the interior of the positive cone is assumed.
\par\smallskip
\begin{enumerate}[resume*=delay_properties]
\item\label{F5}   There exists a pair $(\w_0,\phi_0)\in\W\times \Int\wit K_B$  with bounded semitrajectory $\{(\w_0{\cdot}t, y_t(\w_0,\phi_0))\mid t\ge 0\}$  and such that $y_t(\w_0, \phi_0)\geq_B \psi_0$ for some  $\psi_0\gg_B 0$ and each $t\geq 0$.
\end{enumerate}\par\smallskip
From this  we will show the relatively compactness and uniform stability of each semitrajectory with initial value in $\Int\wit K_B$. Before that, we introduce the part metric on $\Int\wit K_B$ defined~by
\begin{equation}\label{eq:part-metric}
p(\phi,\psi) = \inf\big\{\ln \alpha\mid \alpha\geq 1\;\text{and}\;\, \alpha^{-1}\phi\leq_B\psi\leq_B\alpha\,\phi\big\}
\end{equation}
which turns $(\Int\wit K_B, p)$ into a metric space. The map $p:\Int\wit K_B\times \Int\wit K_B\to \R^+$ is continuous with respect to the product topology induced by the norm $\|\cdot\|_L$, due to the inequality
\[
p(\phi,\psi)\leq \ln \left(1+\frac{\n{\phi-\psi}_L}{R}\right) \quad\text{for each }\phi,\,\psi\in \Int\wit K_B\,,
\]
provided that the closed balls of radius $R>0$ centered at $\phi$ and $\psi$, respectively, are contained in $\wit K_B$ (see Krause and Nussbaum~\cite[Lemma 2.3(i)]{paper:KrNu}). In addition, since $\wit K_B$ is a normal cone, we also have (see~\cite[Lemma 2.3(ii)]{paper:KrNu}):
\begin{equation}\label{ine-norm-p}
\n{\phi-\psi}_L\leq \big(2 \,e^{p(\phi,\psi)}-e^{-p(\phi,\psi)}-1\big)\min(\n{\phi}_L,\n{\psi}_L)\,,\quad \phi,\,\psi\in \Int\wit K_B\,.
\end{equation}

We say that a compact set $M\subset \W\times X_L$ is {\em strongly positive\/} if $M\subset \W\times \Int \wit K_B$, and write $M\gg_B 0$. Note that this is equivalent to the existence of a $\psi_0\gg_B 0$ such that $\psi\geq_B \psi_0$ for all $(\w,\psi)\in M$.
\begin{prop}\label{prop:relcomp}
Assume that conditions \ref{F1}, \ref{F2}, \ref{F4}, and \ref{F5} hold.  Then, for each $\w\in\W$ and $\phi\gg _B 0$, the semiorbit
$\{(\w{\cdot}t, y_t(\w,\phi))\mid t\ge 0\}$ is a relatively compact subset of $\W\times X_L$  and it is  uniformly stable. Moreover, the omega-limit set is strongly positive, i.e.,  $\mathcal{O}(\w,\phi)\gg _B0$, and the restriction of $\tau$ to the compact invariant set $\mathcal{O}(\w,\phi)\subset \W\times X_L$  is continuous.
\end{prop}
\begin{proof}
Once we  check that $\{(\w{\cdot}t, y_t(\w,\phi))\mid t\ge 0\}$  is bounded, it follows easily from condition~\ref{F1}  that it is relatively compact in $\W\times X$. More precisely, given a sequence $\{t_n\}$ there is a subsequence (for simplicity we take the whole sequence) for which $\w{\cdot}t_n\to \w_1$ and $y_{t_n}(\w,\phi)\to \phi_1$  in $X$ as $n\to\infty$. Thus, from
\[
y_{t_n}(\w,\phi)(s)=y_{t_n}(\w,\phi)(-r)+\int_{-r}^s f(\w{\cdot}(t_n+u),y_{t_n+u}(\w,\phi))\,du \,,\quad s\in[-r,0]\,,
\]
we deduce that $\phi_1(s)=\phi_1(-r)+\int_{-r}^s f(\w_1{\cdot}u,(\phi_1)_u)\, du$, i.e., $\phi_1'(s)=f(\w_1{\cdot}{s},(\phi_1)_s)$, which together with
$(y_{t_n}(\w,\phi))'(s)=f(\w{\cdot}(t_n+s),y_{t_n+s}(\w,\phi))\to f(\w_1{\cdot}{s},(\phi_1)_s)$ shows that $(\w{\cdot}t_n,y_{t_n}(\w,\phi))\to (\w_1,\phi_1)$ in $\W\times X_L$, as wanted.
\par
In order to check the boundedness, from~\ref{F5} and the previous argument we deduce that $\{(\w_0{\cdot}t, y_t(\w_0,\phi_0))\mid t\ge 0\}$  is relatively compact in $\W\times X_L$ and, therefore, the omega-limit set $\mathcal{O}(\w_0,\phi_0)$ is a compact and invariant subset of $\W\times\Int\wit K_B$. Moreover, there are $e_1$ and $e_2\in \Int\wit K_B$ such that
\begin{equation}\label{desiw0phi0}
0\ll_B  e_1\leq _B\varphi\leq_Be_2 \quad \text{for each } (\w,\varphi)\in \mathcal{O}(\w_0,\phi_0).
\end{equation}
Then, given $\w\in\W$ and $\phi\gg_B 0$, we can take a pair  $(\w,\varphi)\in  \mathcal{O}(\w_0,\phi_0)$  (thanks to the minimality of $\W$) and a  $\mu>1$ such that $\mu^{-1}\varphi\ll_B\phi\leq_B \mu\,  \varphi$. By monotonicity and sublinearity, $\mu^{-1}y_t(\w,\varphi)\leq _B y_t(\w,\phi)\leq_B \mu\, y_t(\w,\varphi)$ for $t\geq 0$ and the boundedness of the semiorbit of $(\w,\phi)$, both from above and from below,  follows from the invariance of the set $\mathcal{O}(\w_0,\phi_0)$, relation~\eqref{desiw0phi0}, and the semimonotonicity of the norm $\n{\cdot}_L$. Besides, it follows that  $\mathcal{O}(\w,\phi)\gg_B 0$.
\par
As for  the uniform stability of the semiorbit of $(\w,\phi)$, from the uniform continuity of the part metric $p$ over the relatively compact set
$\{ y_s(\w,\phi)\mid s\ge 0\}\times \{ y_s(\w,\phi)\mid s\ge 0\}\subset \Int\wit K_B \times \Int\wit K_B$, given $0<\ep<1$  there is  a $\delta(\ep)\geq 0 $ such that if $\n{y_s(\w,\phi)-\psi}_L<\delta$ for some $s\geq 0$ and some $\psi\in \Int\wit K_B$, then $p(y_s(\w,\phi),\psi) <\ep = \ln (e^{\ep})$, that is,
$e^{-\ep} y_s(\w,\phi)\leq_B\psi\leq_B e^{\ep} y_s(\w,\phi)$, and again  monotonicity and sublinearity  yield
\[
 e^{-\ep} y_{s+t}(\w, \phi)\leq_B y_t(\w{\cdot} s,\psi)\leq_B e^{\ep} y_{s+t}(\w,\phi)\quad \text{for each } t\geq 0\,,
\]
i.e.,  $p(y_t(\w{\cdot}s,\psi),y_{s+t}(\w,\phi))\leq \ln (e^{\ep})=\ep$ for $t\geq 0$, which together with~\eqref{ine-norm-p} provides the uniform stability, as claimed.

Finally, the continuity of $\tau$ restricted to the compact invariant set $\mathcal{O}(\w,\phi)\subset \W\times \Int\wit K_B$ follows from the continuity of $\tau$ in~\eqref{skewC} and  the equivalence of the topologies on $\W\times X$ and  $\W\times X_L$ when restricted to $\mathcal{O}(\w,\phi)$. The proof is finished.
\end{proof}
After this result, a direct application of~\cite[Proposition~3.6(ii)]{noos07JDE} provides a more accurate description of the previous omega-limit sets.
\begin{prop}\label{prop:omega-minimal}
Assume that conditions \ref{F1}, \ref{F2}, \ref{F4}, and \ref{F5} hold. Then, for each $\w\in\W$ and $\phi\gg _B 0$, the omega-limit set $\mathcal{O}(\w,\phi)$  is a uniformly stable and strongly positive minimal set which admits a fiber distal flow extension.
\end{prop}
Even if the semiflow $\tau:\R^+\times \W\times X_L\to \W\times X_L$ is not continuous, we can determine an invariant and locally compact subset $E$ of $\W\times \Int\wit K_B$ formed by entire relatively compact trajectories,  over which the restriction of the skew-product semiflow~\eqref{skewC} is continuous for the norm $\n{\cdot}_L$. Before we state this result, we include a technical lemma.
\begin{lema}\label{lema:tecnico}
Assume that conditions \ref{F1}, \ref{F2}, \ref{F4}, and \ref{F5} hold and let  $(\w_0,\phi_0)$ be the pair in condition \ref{F5}.  Then:
\begin{itemize}
\item[{\rm (i)}] There exists a $\mu>1$ such that
$\mu^{-1}\phi_0\leq_B  y_t(\w,\phi_0)\leq _B \mu\,\phi_0$ for all $\w\in\W$  and $t\geq 0$.
\item[{\rm (ii)}]  If $\varphi\in C_\rho$ for  the  bounded set $C_\rho:=\big\{\varphi\in X_L\mid  \rho^{-1}\phi_0 \leq_B \varphi \leq_B \rho\,\phi_0\big\}$ with $\rho\geq 1$, then $y_t(\w,\varphi) \in  C_{\mu\rho}$ for $\w\in\W$, $t\geq 0$, and $\mathcal{O}(\w,\varphi)\subset \W\times C_{\mu\rho}$.
\end{itemize}
\end{lema}
\begin{proof}
Since the semitrajectory of $(\w_0,\phi_0)$ is bounded and strongly apart from $0$, and $\phi_0\gg_B 0$, we can find a $\mu_0>1$ big enough so that ${\mu_0}^{-1}\phi_0\leq_B y_t(\w_0,\phi_0)\leq_B \mu_0\,\phi_0$ for all $t\geq 0$. Now, given any $\w\in\W$ we can take a pair $(\w,\phi)\in \mathcal{O}(\w_0,\phi_0)$ which also satisfies ${\mu_0}^{-1}\phi_0\leq_B \phi\leq_B  \mu_0\,\phi_0$. Applying monotonicity,  sublinearity, and the fact that $\tau_t(\w,\phi)$ remains in $\mathcal{O}(\w_0,\phi_0)$ for all $t\geq 0$, we have that  ${\mu_0}^{-1}y_t(\w,\phi_0)\leq_B y_t(\w,\phi)\leq_B  \mu_0\,\phi_0$  and also ${\mu_0}^{-1}\phi_0\leq_B y_t(\w,\phi)\leq_B  \mu_0\,y_t(\w,\phi_0)$  for all $t\geq 0$.  As a consequence, the proof of (i) is finished by taking  $\mu=\mu_0^2$.

Item (ii) follows straighforward by applying monotonicity, sublinearity, and (i). The proof is finished.
\end{proof}
We keep the terminology just introduced for the sets $C_\rho\subset X_L$ for $\rho\geq 1$.
\begin{theorem}\label{teor:K}
Assume that conditions \ref{F1}, \ref{F2}, \ref{F4}, and \ref{F5} hold  and consider the set $
E=\bigcup_{\substack{\w\in\W\\ \phi\gg_B 0}}\mathcal{O}(\w,\phi)$.
Then:
\begin{itemize}
\item[\rm{(i)}] The set $E\subset \W\times \Int\wit K_B\subset \W\times X_L$ is invariant and locally compact.
\item[\rm{(ii)}]  $\tau\colon \R^+\times E\longrightarrow  E$  is continuous for the topology of $\W\times X_L$ on $E$.
\item[\rm{(iii)}] For each $j\geq 1$,
\begin{align*}
&\lim_{t\to\infty}{\rm dist}(y_t(\w,C_j),E(\w{\cdot}t))=0\quad\text{and}\\
&\lim_{t\to \infty}{\rm dist}(y_t(\w{\cdot}(-t),C_j),E(\w))=0
\end{align*}
uniformly for $\w\in\W$, where as usual $E(\w)$ denotes the $\w$-section of $E$.
\end{itemize}
\end{theorem}
\begin{proof}
First of all notice that $E\subset \W\times \Int\wit K_B\subset \W\times X_L$ by Proposition~\ref{prop:relcomp} and $E$ is invariant, that is, $\tau_t(E)=E$ for all $t\geq 0$. To see that $E$ is locally compact,  let us consider the family of open subsets in $E$ given by
\begin{equation}\label{eq:Vj}
V_j:=\left\{(\w,\varphi)\in E \,\left|\; \frac{1}{j}\,\phi_0 \ll_B \varphi \ll_B j\,\phi_0 \right.\right\} \quad\text{for each}\;j\geq 1\,,
\end{equation}
 where $(\w_0,\phi_0)$ is the pair in condition \ref{F5}.
Note that, given any $(\wit\w,\wit\varphi)\in E$, there is a sufficiently big $j$ such that $(\wit\w,\wit\varphi)\in V_j$. Thus, if we prove that the sets $V_j$ are relatively compact in $\W\times X_L$, we are done.

With this purpose, fix a $j\geq 1$ and a sequence $\{(\w_n,\varphi_n)\} \subset V_j$ and let us check that it has a convergent subsequence. Note that, for each $n\geq 1$, $(\w_n,\varphi_n)\in \mathcal{O}(\wit\w_n,\wit\varphi_n)$ for some $\wit\w_n\in\W$ and $\wit\varphi_n\gg_B 0$. By Proposition~\ref{prop:omega-minimal}, the set $\mathcal{O}(\wit\w_n,\wit\varphi_n)$ admits a flow extension, that is, there exists a unique backward orbit for each of its points. Then, there exists a unique $\psi_n\in X_L$ such that $\tau(r,\w_n{\cdot}(-r), \psi_n)=(\w_n,\varphi_n)$. If the sequence $\{\psi_n\}$ were relatively compact in $X$ for the sup-norm, then for a subsequence $\psi_{n_k}\to \psi_0\in X$ as $k\to\infty$. Since $\W$ is compact, we can assume w.l.o.g.~that $\w_{n_k}\to \w\in \W$ as $k\to\infty$. Therefore, by the continuity of the map $\tau_r:\W\times X\to \W\times X_L$ (see~\eqref{XXL}),  $y_r(\w_{n_k}{\cdot}(-r), \psi_{n_k})= \varphi_{n_k} \to y_r(\w{\cdot}(-r), \psi_0)$ in $X_L$ as $k\to\infty$ and we would have found a convergent subsequence $\{(\w_{n_k},\varphi_{n_k})\}$, as wanted.

So, to finish the proof of (i) we need to see that $\{\psi_n\}$ is relatively compact in $X$ for the sup-norm.
The main idea is to choose appropriate pairs $(\w_n^*,\varphi_n^*)\in \W\times \Int\wit K_B$ in such a way that  $(\w_n,\varphi_n)\in \mathcal{O}(\w_n^*,\varphi_n^*)$  for each $n\geq 1$ and the set
\begin{equation*}
H^*=\{y_{t+r}(\w_n^*,\varphi_n^*)\mid n\geq 1,\,t\geq 0\}\subset X
\end{equation*}
is relatively compact. Note that then, also $(\w_n{\cdot}(-r), \psi_n) \in \mathcal{O}(\w_n^*,\varphi_n^*)$ for each $n\geq 1$, and thus $\{\psi_n\}\subset \cls H^*$, the closure of $H^*$ under the sup-norm, a compact set.

By the definition of $V_j$, for each $n\geq 1$, we have that $j^{-1}\phi_0 \ll_B \varphi_n \ll_B j\,\phi_0$
and, since $(\w_n,\varphi_n)\in \mathcal{O}(\wit\w_n,\wit\varphi_n)$, for each $n\geq 1$ there exists a time $s_n$ big enough so that
$j^{-1}\phi_0 \ll_B y_{s_n}(\wit\w_n,\wit\varphi_n) \ll_B j\,\phi_0$. Let $(\w_n^*,\varphi_n^*)=\tau(s_n,\wit\w_n,\wit\varphi_n)$. It is clear that $(\w_n,\varphi_n)\in \mathcal{O}(\w_n^*,\varphi_n^*)$ for each $n\geq 1$ and $\varphi_n^* \in C_j$.  By Lemma~\ref{lema:tecnico},  $y_t(\w_n^*,\varphi_n^*) \in C_{\mu j}$ for all $t\geq 0$ and $n\geq 1$, so that in particular $y_t(\w_n^*,\varphi_n^*)$ is in the interval for the usual order in $X$, $I=[j^{-1}\mu^{-1}\,\phi_0,j\,\mu\,\phi_0]$.
Then, by~\eqref{semicocycle}, we can rewrite
\[
H^*=\{y_r(\w_n^*{\cdot}t,y_t(\w_n^*,\varphi_n^*))\mid n\geq 1,\,t\geq 0\}\subset y_r\big(\W\times I\big)\,,
\]
which is relatively compact in $X$ because the map $y_r$ is compact. Therefore, $H^*$ is relatively compact too and the proof of (i) is finished.

Most of the arguments for the proof of (ii) are identical to the previous ones, so that we just give a sketch of the proof.
To check the joint continuity of $\tau$ on $\R^+\times E$, we take $\{(t_n,\w_n,\varphi_n)\}\subset \R^+\times E$ converging to $(\bar t,\bar\w,\bar\varphi)\in \R^+\times E$ as $n\to\infty$ with the norm $\|\cdot\|_L$ in $X_L$. We already know that $\lim_{n\to\infty}\w_n{\cdot}t_n=\bar\w{\cdot}\bar t$.
For the second component of $\tau$, note that, for each $n\geq 1$, $(\w_n,\varphi_n)\in \mathcal{O}(\wit\w_n,\wit\varphi_n)$ for some $\wit\w_n\in\W$ and $\wit\varphi_n\gg_B 0$. Again thanks to~\eqref{semicocycle} we can write
\begin{equation}\label{eq:y_tn}
y_{t_n}(\w_n,\varphi_n)= y_{t_n+r}(\w_n{\cdot}(-r), \psi_n)\,,
\end{equation}
where $\psi_n\in X_L$  satisfies $y_r(\w_n{\cdot}(-r), \psi_n)=\varphi_n$. If the sequence $\{\psi_n\}$ were relatively compact in $X$ for the sup-norm, then $\psi_n\to \psi_0\in X$ as $n\to\infty$, up to taking a subsequence. By the continuity of $\tau: [r,\infty)\times\W\times X\to\W\times X_L$,   we would first have  that $y_r(\w_n{\cdot}(-r), \psi_n)=\varphi_n\to y_r(\bar\w{\cdot}(-r), \psi_0)=\bar\varphi$  in $X_L$ as $n\to\infty$, and then, taking limits in~\eqref{eq:y_tn},   $y_{t_n}(\w_n,\varphi_n)\to y_{\bar t+r}(\bar \w{\cdot}(-r), \psi_0)=y_{\bar t}(\bar \w, y_r(\bar\w{\cdot}(-r), \psi_0))=y_{\bar t}(\bar \w,\bar\varphi)$ also in $X_L$ as $n\to\infty$, and we would be done.

So, it remains to check that $\{\psi_n\}$ is relatively compact in $X$, and this is done just as in the proof of (i). We just point out that this time we can take a big $j\geq 1$ and a ball $D$ centered at $\bar\varphi\gg_B 0$ so that  ${j}^{-1}\phi_0\leq_B \phi \leq_B j\,\phi_0$ for all  $\phi\in D$. Since $\varphi_n\to\bar\varphi$ as $n\to\infty$ and $(\w_n,\varphi_n)\in \mathcal{O}(\wit\w_n,\wit\varphi_n)$, for each $n\geq 1$ we can find a time  $s_n$ big enough so that $y_{s_n}(\wit\w_n,\wit\varphi_n)\in D$.  From here the proof is finished just as before.

In the proof of (iii) we will use some properties of the Hausdorff semidistance \eqref{semidist}. Namely, ${\rm dist}(Y_1,Y_2)={\rm dist}(Y_1,\overline{Y}_{\! 2})$; if ${\rm dist}(Y_1,Y_2)=0$, then $Y_1\subseteq \overline{Y}_{\! 2}$; if $Y_1\subset Z_1$, then ${\rm dist}(Y_1,Y_2)\leq {\rm dist}(Z_1,Y_2)$; and if $Z_2\subset Y_2$, then ${\rm dist}(Y_1,Y_2)\leq {\rm dist}(Y_1,Z_2)$, for subsets $Y_1, Y_2, Z_1, Z_2$ of a metric space $Y$, denoting by  $\overline{Y}_{\! 2}$ the closure of the set, for the sake of notation. Finally, recall that when $Y_1$ is a singleton, the Hausdorff semidistance is just the usual distance between a point and a set, ${\rm d}(y_1,Y_2)$.

The proof is organised into a series of claims. Let us fix a $j\geq 1$ and let $\mu>1$ be the one in Lemma~\ref{lema:tecnico}.
\par\smallskip
\textbf{Claim 1.\/} Let $j_0\geq \mu\,j$ and consider the set $E_{j_0}:=\bigcup_{(\w,\varphi)\in\W\times C_{j_0}}\mathcal{O}(\w,\varphi)\subset E$.
Then, $E_{j_0}$ is invariant and relatively compact in $\W\times X_L$, the section map $\w\in\W\mapsto \overline E_{j_0}(\w)$ is continuous at every point of $\W$, and for each fixed $\wit\w\in\W$, $\lim_{t\to\infty}{\rm dist}(y_t(\wit\w,C_j),E_{j_0}(\wit\w{\cdot}t))=0$.
\par\smallskip
First of all, note that by taking $E_{j_0}$ we are reducing the zone of the set $E$ to which $y_t(\w,C_j)$ is going to approach, with the aim that it is  relatively compact and uniformly stable. Recall that if $\varphi\in C_j$, then $\mathcal{O}(\w,\varphi)\subset \W\times C_{\mu j}$ for $\w\in\W$, and $C_{\mu j}\subseteq  C_{j_0}$. Iterating this fact,  taking $j_1> \mu\,j_0$, $E_{j_0}\subset V_{j_1}$ for the relatively compact set $V_{j_1}$ defined  in \eqref{eq:Vj}, so that $E_{j_0}$ is relatively compact too.
Besides, arguing as in the proof of Proposition~\ref{prop:relcomp} we get the uniform stability of this set, that is, given any $\varepsilon>0$ there is a $\delta>0$ such that if $(\w,\phi)\in E_{j_0}$ and $\psi\in \Int\wit K_B$ satisfy $\|\phi-\psi\|_L<\delta$, then $\|y_t(\w,\phi)-y_t(\w,\psi)\|_L<\varepsilon$ for $t\geq 0$. Note also that $E_{j_0}$ admits a flow extension, by Proposition~\ref{prop:omega-minimal}.  Then, we can apply \cite[Theorem~3.3]{noos07JDE} to the invariant compact set $\overline E_{j_0}\subset \W\times X_L$, whence the section map $\w\in\W\mapsto \overline E_{j_0}(\w)$ is continuous at every point of $\W$ for the Hausdorff metric, which implies its continuity also for the Hausdorff semidistance we are using.

At this point,  fix an $\wit\w\in\W$. To see that  $\lim_{t\to \infty} {\rm dist}(y_t(\wit \w,C_j),E_{j_0}(\wit\w{\cdot}t))=0$, let us argue by contradiction and assume that there exists an $\varepsilon_0>0$ and there are sequences $\{t_n\}\uparrow \infty$ and  $\{\varphi_n\}\subset C_j$ such that
\begin{equation}\label{eq:e0}
{\rm d}(y_{t_n}(\wit\w,\varphi_n),  E_{j_0}(\wit\w{\cdot}t_n))\geq \varepsilon_0 \quad\text{for}\; n\geq 1\,.
\end{equation}
Now, as in the proof of Theorem~\ref{teor:atractores}, the set $H=\cls\{y_r(\w,\varphi)\mid (\w,\varphi)\in \W\times C_{j_0}\}$ is compact in $X$ and
$\wit H=\{y_r(\w,\varphi)\mid (\w,\varphi)\in \W\times H\}$ is compact in $X_L$. Besides,  by using the semicocycle property~\eqref{semicocycle},
\begin{equation}\label{eq:wit H}
y_t(\w,\varphi)\in \wit H \quad \text{whenever }\;\w\in \W\,,\;\varphi\in C_j\;\text{ and }\; t\geq 2\,r\,.
\end{equation}
Then, supposing that $t_n\geq 2r$ for $n\geq 1$, we have that $y_{t_n}(\wit\w,\varphi_n)\in \wit H$ for $n\geq 1$ and we can assume w.l.o.g.~that $(\wit\w{\cdot}t_n,y_{t_n}(\wit\w,\varphi_n))\to (\w^*,\varphi^*)$ as $n\to\infty$. Due to the continuity of the section map $\w\in\W\mapsto \overline E_{j_0}(\w)$ at $\w^*$, taking limits in~\eqref{eq:e0} we conclude that ${\rm d}(\varphi^*,\overline E_{j_0}(\w^*))\geq \varepsilon_0$.

On the other hand, since also $\{y_{2r}(\wit\w,\varphi_n)\}\subset \wit H$, once more we can assume w.l.o.g. that $y_{2r}(\wit\w,\varphi_n)\to \wit \varphi\in \wit H \cap  C_{j_0}$ as $n\to\infty$.
Arguing again  as in the proof of Proposition~\ref{prop:relcomp}, we obtain the uniform stability of the compact set $\wit H\subset \Int \wit K_B$, that is, given any $\varepsilon>0$ there is a $\delta>0$ such that, if $\phi\in \wit H$ and $\psi\in \Int\wit K_B$ satisfy $\|\phi-\psi\|_L<\delta$, then $\|y_t(\w,\phi)-y_t(\w,\psi)\|_L<\varepsilon$ for all $\w\in\W$ and $t\geq 0$. Then, given any $\varepsilon>0$ there is an $n_{\varepsilon}$ such that  $\|y_{t_n-2r}(\wit\w{\cdot}2r,y_{2r}(\wit\w,\varphi_n))-y_{t_n-2r}(\wit\w{\cdot}2r,\wit\varphi)\|_L<\varepsilon$ for  $n\geq n_{\varepsilon}$, that is,
$\|y_{t_n}(\wit\w,\varphi_n)-y_{t_n-2r}(\wit\w{\cdot}2r,\wit\varphi)\|_L<\varepsilon$ for $n\geq n_{\varepsilon}$. It follows that $\lim_{n\to\infty} \tau_{t_n-2r}(\wit\w{\cdot}2r, \wit\varphi)=(\w^*,\varphi^*)$, that is, $(\w^*,\varphi^*)\in \mathcal{O}(\wit\w{\cdot}2r,\wit\varphi)\subset E_{j_0}$ by the construction. Therefore, ${\rm d}(\varphi^*,E_{j_0}(\w^*))=0$, a contradiction.

\par\smallskip
\textbf{Claim 2.\/} The {\it omega-limit set of $\W\times C_j$}, $\mathcal{O}_j:=\{(\w,\varphi)\in \W\times X_L \mid (\w,\varphi)=\lim_{n\to\infty}  (\w_n{\cdot}t_n,y_{t_n}(\w_n,\varphi_n))$ for some $\{\w_n\}\subset \W$, $\{\varphi_n\}\subset C_j$, $\{t_n\}\uparrow \infty\}$,
is an invariant compact set and the section map $\w\mapsto \mathcal{O}_j(\w)$ is continuous at every $\w\in\W$.
\par\smallskip
By the construction, $\mathcal{O}_j$ is a closed set. It is positively invariant because of the semiflow property and the continuity of  the map $\tau_t:\W\times X_L\to \W\times X_L$ for each fixed $t\geq 0$. Besides, according to~\eqref{eq:wit H},  $\mathcal{O}_j\subset \W\times \wit H$, so that also $\mathcal{O}_j$ is compact and uniformly stable.  Last but not least, there are backward extensions inside $\mathcal{O}_j$. With all these properties, we can apply~\cite[Theorem~3.4]{noos07JDE} to get the continuity of the section map.
\par\smallskip
\textbf{Claim 3.\/} $\lim_{t\to\infty}{\rm dist}(y_t(\w,C_j),\mathcal{O}_j(\w{\cdot}t))=0$ uniformly for $\w\in \W$.
\par\smallskip
To see it,  argue by contradiction and assume that there exist an $\varepsilon_0>0$ and sequences $\{t_n\}\uparrow \infty$, $\{\w_n\}\subset \W$ and $\{\varphi_n\}\subset C_j$ such that
${\rm d}(y_{t_n}(\w_n,\varphi_n), \mathcal{O}_j(\w_n{\cdot}t_n))\geq \varepsilon_0$ for $n\geq 1$.
Once more, by~\eqref{eq:wit H} we can assume w.l.o.g.~that  $(\w_n{\cdot}t_n,y_{t_n}(\w_n,\varphi_n))\to (\w^*,\varphi^*)\in \mathcal{O}_j$. Thus, on the one hand ${\rm d}(\varphi^*, \mathcal{O}_j(\w^*))=0$, but on the other hand, by the continuity of the section map, ${\rm d}(\varphi^*, \mathcal{O}_j(\w^*))\geq \varepsilon_0$, which is a contradiction.

\par\smallskip
\textbf{Claim 4.\/} Let $j_0\geq \mu\,j$ and $j_1\geq \mu\,j_0$. Then, $\mathcal{O}_j(\w)\subseteq \overline E_{j_1}(\w)$ for all $\w\in\W$.
\par\smallskip
Note that by its definition, $\mathcal{O}_j\subset \W\times C_{j_0}$. Since $\mathcal{O}_j$ is an invariant set, for a fixed $\wit\w\in\W$,   $\mathcal{O}_j(\wit\w{\cdot}t) \subset y_t(\wit\w,C_{j_0})$ for all $t\geq 0$. Therefore, ${\rm dist}(\mathcal{O}_j(\wit\w{\cdot}t),E_{j_1}(\wit\w{\cdot}t))\leq
{\rm dist}(y_t(\wit\w,C_{j_0}),E_{j_1}(\wit\w{\cdot}t))\to 0$ as $t\to\infty$ by Claim~1 applied to $\wit\w$, $j_0$ and $j_1$.

From here, by the minimal character of $\W$, for each $\w\in\W$ we can find a sequence $\{t_n\}\uparrow \infty$ such that $\wit\w{\cdot}t_n\to \w$. By the  continuity of the section maps expressed in Claims~1 and~2, ${\rm dist_{\mathcal H}}(\mathcal{O}_j(\w), \overline E_{j_1}(\w))=\lim_{n\to\infty}{\rm dist_{\mathcal H}}(\mathcal{O}_j(\wit\w{\cdot}t_n),\overline E_{j_1}(\wit\w{\cdot}t_n))=0$. As a consequence, $\mathcal{O}_j(\w)\subseteq \overline E_{j_1}(\w)$ for all $\w\in\W$.

\par\smallskip
\textbf{Claim 5.\/} $\lim_{t\to\infty}{\rm dist}(y_t(\w,C_j),E(\w{\cdot}t))=0$ uniformly for $\w\in \W$.
\par\smallskip
For each $\w\in\W$, $E_{j_1}(\w) \subset E(\w)$ and, by Claim~4, $\mathcal{O}_j(\w)\subset \overline E_{j_1}(\w)$. Then, for all $t\geq 0$, ${\rm dist}(y_t(\w,C_j),E(\w{\cdot}t))\leq  {\rm dist}(y_t( \w,C_j),E_{j_1}(\w{\cdot}t))\leq {\rm dist}(y_t(\w,C_j),\mathcal{O}_j(\w{\cdot}t))$, and this goes to $0$ as $t$ goes to $\infty$ uniformly for $\w\in\W$, as stated in Claim~3.
\par\smallskip
Finally, note that the second limit  in (iii) follows straightforward from the first one. The proof is finished.
\end{proof}
\begin{remark}
Note that the set $E$ in Theorem~\ref{teor:K} is formed by bounded and entire orbits and, when there exists a global attractor in $\W\times \Int\wit K_B$, it coincides with this set. This happens, for instance, under the conditions in Theorem~\ref{teor:atractores}(ii). In any case, as stated in (iii), the set $E$ contains the essential ingredients of the pullback and forwards dynamics of the strongly positive semitrajectories.
\end{remark}
Finally, the existence of a point $\w_1\in\W$ and a time $t_1>0$ of strong sublinearity for the trajectories,  which  in applications will be checked in terms of one of the FDEs of the family, allows us to characterise the long-term behaviour of the trajectories.
\par\smallskip
\begin{enumerate}[resume*=delay_properties]
\item\label{F6}   There exist a point $\w_1\in\W$ and a time  $t_1>0$ such that
\[
y_{t_1}(\w_1,\lambda\,\psi)\gg_B \lambda\, y_{t_1}(\w_1,\psi) \quad\text{whenever }\psi\gg_B 0 \text{ and } \lambda \in(0,1) \,.
\]
\end{enumerate}
\begin{theorem}\label{th:main}
Assume that conditions \ref{F1}, \ref{F2}, and \ref{F4}--\ref{F6} hold. Then, there exists a unique strongly positive minimal set $M\gg_B 0$ which is a copy of the base, that is, $M=\{(\w,b(\w))\mid \w\in \W\}$ for a continuous map $b:\Omega\to \Int\wit K_B$, and
\begin{equation}\label{convergence-to-M}
\lim_{t\to\infty}\|y_t(\w,\phi)-b(\w{\cdot}t)\|_L=0 \quad \text{whenever } \w \in \W \text{ and }\phi\gg_B 0\,.
\end{equation}
\end{theorem}
\begin{proof}
First notice that there is at least one strongly positive minimal set because, if we take the point $(\w_0,\phi_0)$ in condition~\ref{F5},  Proposition~\ref{prop:omega-minimal} asserts that $M=\mathcal{O}(\w_0,\phi_0)$ is a strongly positive minimal set, $M\gg_B 0$. Assume on the contrary that there are two distinct minimal sets $M_1\gg_B 0$ and $M_2\gg_B 0$, and consider
$(\w,\psi_1)\in M_1$ and  $(\w,\psi_2)\in M_2$.  Let $\w_1\in\W$ and $t_1>0$ be the point and time of  condition~\ref{F6}. We take a sequence $\{s_n\}\downarrow -\infty$ with $s_{n-1}-s_n>t_1$ for each $n\geq 2$, and such that
\begin{equation}\label{eq:a-inf}
 \lim_{n\to\infty}(\w{\cdot}s_n, y_{s_n}(\w,\psi_i))=(\w_1,\psi^*_i)\in M_i \;\text{ for }\; i=1,2\,.
\end{equation}
Note that in particular $\psi^*_1\not=\psi^*_2$, since distinct minimal sets have an empty intersection. Now, from~\ref{F2}, \ref{F4}, and~\ref{F6}, i.e., monotonicity, sublinearity, and strong sublinearity at time $t_1$ for  $\w_1$, it is well-known (see,  e.g., Chueshov~\cite[Lemma~4.2.1]{book:chue})  that the part  metric $p$ defined in~\eqref{eq:part-metric} is decreasing along the trajectories, that is,
\[
p(\varphi_1,\varphi_2)\geq p(y_t(\w,\varphi_1),y_t(\w,\varphi_2)) \quad \text{whenever } \w\in\W,\;\varphi_1,\varphi_2\in  \Int\wit K_B \text{ and }t\geq 0\,,
\]
and strictly decreasing at $t_1$ for the point  $\w_1$, that is,
\[
p(\varphi_1,\varphi_2)>p(y_{t_1}(\w_1,\varphi_1),y_{t_1}(\w_1,\varphi_2))\quad \text{whenever } \varphi_1,\varphi_2\in  \Int\wit K_B,\;\varphi_1\not=\varphi_2\,.
\]
This fact, together with~\eqref{eq:a-inf}, the continuity of the part metric,  and the inequalities $s_{n-1}> t_1+s_n$ for $n\geq 2$, yield
\begin{align*}
p(\psi^*_1,\psi^*_2)& >p(y_{t_1}(\w_1,\psi^*_1),y_{t_1}(\w_1,\psi^*_2)) =\lim_{n\to\infty} p(y_{t_1+s_n}(\w,\psi_1),y_{t_1+s_n}(\w,\psi_2))\\
                    &  \geq \lim_{n\to\infty} p(y_{s_{n-1}}(\w,\psi_1),y_{s_{n-1}}(\w,\psi_2))= p(\psi^*_1,\psi^*_2)\,,
\end{align*}
a contradiction. Thus, there is a unique strongly positive minimal set $M\gg_B 0$ and now we check that it is a copy of the base. Once more we argue by contradiction and assume that for certain $\w\in \W$ there are two distinct pairs $(\w,\psi_1), (\w,\psi_2)\in M$. Then, we get a contradiction arguing as before. Just note that this time we get that $\psi^*_1\not=\psi^*_2$ because $M$ has a fiber distal flow extension (see Proposition~\ref{prop:omega-minimal}). Consequently $M$ is a copy of the base, $b$ is continuous,   and~\eqref{convergence-to-M} follows from $\mathcal{O}(\w,\phi)=M$ for all $\w\in\W$ and $\phi\gg_B 0$. The proof is finished.
\end{proof}
We remark that the implications of Theorem~\ref{th:main} in terms of attraction are strong: the induced semiflow has a simple global attractor in $\W\times \Int \wit K_B$, and for all the processes the pullback attractor is a forwards attractor too.
\section{An application to Nicholson systems}\label{sec-nicholson}
In this section we want to find some applications of the theory to delay systems in real life processes. We have focused on Nicholson systems, modelling the behaviour of a biological species on an heterogeneous environment,  giving rise to patches or compartments in the model, so that the distribution of the population is influenced by the migrations among patches and the growth of the population on each patch, which depends on the local resources, among other conditions.
For the sake of simplicity, and motivated by the recent interest in almost periodic Nicholson systems, we assume an almost periodic variation of  the coefficients, but all the results might be stated in a more general framework. Namely, it suffices to assume a recurrent behaviour of the time-varying coefficients, so that the hull of the system is minimal. It is noteworthy that,  in that case, the characterisation of  the persistence properties for almost periodic Nicholson systems given in the papers by Obaya and Sanz~\cite{obsa16,obsa18} takes a slightly more complicated form.
\par
Some history on the Nicholson models can be found in many papers, for instance in~\cite{obsa18}. The first model was presented by Gurney et al.~\cite{gubl}, who proposed the scalar delay equation
\begin{equation*}
x'(t)=-\mu\,x(t)+p\,x(t-r)\,e^{-\gamma\,x(t-r)}\,,
\end{equation*}
which was called Nicholson's blowflies equation, as it suited reasonably well the experimental data on the behaviour of an Australian sheep-blowfly obtained by Nicholson. The coefficients in the equation are positive constants with a biological interpretation, and in particular the delay  $r$ stands for the maturation time of the species. When the  equation is rescaled into $x'(t)=-a\,x(t)+b\,x(t-1)\,e^{-x(t-1)}$, it is well-known that,  if $1<b/a\leq e$, then the nontrivial equilibrium attracts every other nontrivial positive solution (see~\cite[Theorem~6.5.1]{smit}). This behaviour was extended later in~\cite{faro} under the condition $1<b/a\leq e^2$. Smith~\cite{smit} introduced the exponential ordering in the study of this equation, as it permits to add a new zone for the variation of the parameters, namely, $b/a>e$ and $b<e^{1-a}$, still guaranteeing that the nontrivial equilibrium of this scalar equation attracts every other nontrivial positive solution.
\par
Our purpose is to apply the exponential ordering in the search for new conditions to ensure the existence of a unique attracting positive almost periodic solution of an almost periodic Nicholson system.
Special attention has been paid in the literature to the periodic case, with quite successful results (see~\cite{faria17} and the references therein), but actually reality is better modelled by almost periodic variations rather than periodic ones.  For convenience,  recall that a continuous map $f:\R\to\R$ is almost periodic if for every $\varepsilon>0$ the set of so-called $\varepsilon$-{\em periods\/} of $f$, $\{s\in\R\mid |f(t+s)-f(t)|<\varepsilon \;\text{for all}\; t\in\R\}$, is relatively dense.
\par
More precisely, we consider an $m$-dimensional system of delay FDEs with patch structure ($m$ patches) and a nonlinear term of Nicholson type, which reflects an almost periodic temporal variation in the environment,
\begin{equation}\label{nicholson delay}
y_i'(t)=-\wit d_i(t)\,y_i(t) +\des \sum_{j=1}^m \wit a_{ij}(t)\,y_j(t) + \wit\beta_{i}(t)\,y_i(t-r_i)\,e^{-\wit c_i(t)\,y_i(t-r_i)}\,,\quad t\geq 0\,
\end{equation}
for $i=1,\ldots,m$. Here $y_i(t)$ denotes the density of the population on patch $i$ at time $t\geq 0$ and
$r_i>0$ is the maturation time on that patch. We make the following assumptions on the coefficient functions:
\par\smallskip
\begin{enumerate}[label=\upshape(\text{a$\arabic*$}),series=nicholson_properties,leftmargin=27pt]
\item\label{a1} $\;\wit d_i(t)$, $\wit a_{ij}(t)$, $\wit c_i(t)$ and $\wit \beta_{i}(t)$ are almost periodic maps on $\R$;
\item\label{a2} $\;\wit d_i(t)\geq d_0>0$ for each $t\in\R$ and $i\in\{1,\ldots,m\}$;
\item\label{a3} $\;\wit a_{ij}(t)$ are all nonnegative maps and $\wit a_{ii}$ is  identically null;
\item\label{a4} $\;\wit\beta_{i}(t)>0$ for each $t\in \R$ and $i\in\{1,\ldots,m\}$;
\item\label{a5} $\;\wit c_i(t)\geq c_0>0$ for each $t\in\R$ and $i\in\{1,\ldots,m\}$;
\item\label{a6} $\;\wit d_i(t)-\sum_{j=1}^m \wit a_{ji}(t)>0$  for each $t\in \R$ and $i\in\{1,\ldots,m\}\,$.
\end{enumerate}
\par\smallskip\noindent
The coefficient $\wit a_{ij}(t)$ stands for the migration rate of the population moving from patch $j$ to patch
$i$ at time $t\geq 0$. As for the birth function, it is given by the delay Nicholson term. Condition~\ref{a5} is technical and implies the uniform boundedness of the terms $y\,e^{-\wit c_i(t)\,y}$ for $y\geq 0$, $t\in\R$, $1\leq i\leq m$.  Finally, the decreasing rate on patch $i$, given by $\wit d_i(t)$, includes the mortality rate as well as the migrations coming out of patch $i$, so that condition~\ref{a6} makes sense, saying that the mortality rate is positive at every time.
\par
We note that we need coefficients defined on  $\R$ to build the {\em hull\/} of the Nicholson system. The reason for introducing the family of systems over the hull is that it allows the use of techniques of non-autonomous dynamical systems.  We briefly explain the procedure. Take $X=C([-r_1,0])\times \ldots \times C([-r_m,0])$ with the usual cone of positive elements, denoted by $X_+$, and the sup-norm. Then, $X$ is a strongly ordered Banach space. Write~\eqref{nicholson delay} as $y_i'(t)=f_i(t,y_t)$, $1\leq i\leq m$, for the maps $f_i:\R\times X\to \R$,
\begin{equation}\label{eq:f}
f_i(t,\phi)=- \wit d_i(t)\,\phi_i(0) +\des \sum_{j=1}^m  \wit a_{ij}(t)\,\phi_j(0) + \wit \beta_{i}(t)\,\phi_i(-r_i)\,e^{-\wit c_i(t)\,\phi_i(-r_i)}\,.
\end{equation}
Consider the map $l:\R\to\R^N$ given by all the almost periodic coefficients $l(t)=(\wit d_i(t),\wit a_{ij}(t),\wit \beta_i(t),\wit c_i(t))$ and let $\W$ be its hull, that is, the closure of the time-translates of $l$ for the compact-open topology. Then, $\W$ is a compact metric space thanks to the boundedness and uniform continuity of almost periodic maps. Besides, the shift map $\sigma:\R\times \W\to\W$, $(t,\w)\mapsto \w{\cdot}t$, with $(\w{\cdot}t)(s)= \w(t+s)$, $s\in\R$, defines an almost periodic and minimal flow. By considering the continuous nonnegative maps $d_i,\,a_{ij},\,\beta_{i},\,c_i:\W\to\R$ such that $(d_i(\w), a_{ij}(\w), \beta_i(\w), c_i(\w))=\w(0)$, the initial system is included in the family of systems over the hull, which can be written for each $\w\in\W$~as
\begin{equation}\label{nicholson delay hull}
y'_i(t)=- d_i(\w{\cdot}t)\,y_i(t) +\des \sum_{j=1}^m  a_{ij}(\w{\cdot}t)\,y_j(t) + \beta_{i}(\w{\cdot}t)\,y_i(t-r_i)\,e^{-c_i(\w{\cdot}t)\,y_i(t-r_i)}\,
\end{equation}
for $i=1,\ldots,m$.   For each $\w \in \W$ and $\varphi\in X$, the solution of~\eqref{nicholson delay hull} with initial value $\varphi$ is denoted by $y(t,\w,\varphi)$.  Solutions  induce a skew-product semiflow $\tau:\R_+\times \W\times X\to \W\times X$, $(t,\w,\varphi)\mapsto (\w{\cdot}t,y_t(\w,\varphi))$ (in principle only locally defined) which has a trivial minimal set $\W\times \{0\}$, as the null map is a solution of all the systems over the hull.

Note that this family of systems does not satisfy the standard  quasimonotone condition~\ref{Q} in FDEs. In any case, the set $\W\times X_+$ is invariant for the dynamics, that is, the solutions of~\eqref{nicholson delay hull} starting inside the positive cone remain inside the positive cone while defined: just apply the criterion given in~\cite[Theorem~5.2.1]{smit}. We can also assert that, if $\varphi\geq 0$ with $\varphi(0)\gg 0$, then $y(t,\w,\varphi)\gg 0$ for all $t\geq 0$. In order to check it, just  compare the solutions of the Nicholson systems with those of the cooperative family of ODEs given for each  $\w\in\W$ by
\begin{equation*}
y'_i(t)=- d_i(\w{\cdot}t)\,y_i(t) +\des \sum_{j=1}^m  a_{ij}(\w{\cdot}t)\,y_j(t)\,,\quad 1\leq i\leq m \,.
\end{equation*}
\par
As stated in~\cite[Theorem~3.3]{obsa18}, all the solutions of~\eqref{nicholson delay hull} are ultimately bounded, that is, there exists a constant $r_0>0$ such that, for each $\w\in \W$ and  $\varphi\in X_+$,  every component of the vector solution satisfies $0\leq y_i(t,\w,\varphi)\leq r_0$ from one time on. As a consequence, the induced semiflow is globally defined on $\W\times X_+$.
\par
Hereafter we focus on situations in which the population persists. We give the definition of persistence for~\eqref{nicholson delay}, meaning that, if at the initial time $t=0$ there are some individuals on every patch,  in the long run  the population will surpass a positive lower bound on all the patches.  We keep the terminology used in~\cite{obsa18}.
\begin{definition}\label{defi persistence Nicholson}
The Nicholson system~\eqref{nicholson delay} is {\em uniformly persistent at $0$}  if there exists an $M>0$ such that for every initial map $\varphi\geq 0$ with  $\varphi(0)\gg 0$ there exists a time $t_0=t_0(\varphi)$ such that
     \[ y_i(t,\varphi)\geq M \quad \text{for all }\;t\geq t_0 \;\text{ and  }\; i=1,\ldots,m \,.\]
\end{definition}
As shown in~\cite[Theorem~3.4]{obsa18}, this dynamical property for the system is equivalent to the uniform persistence of the whole family
\eqref{nicholson delay hull}, or in other words, to the uniform persistence of the induced skew-product semiflow $\tau$ in the interior of the positive cone, according to Definition~\ref{defi persistence hull}(i). A complete characterization of this property in terms of a few Lyapunov exponents (which can be numerically calculated) can be found in~\cite[Theorem~3.5]{obsa18}.

We remark that the equivalence between the uniform persistence of the initial system and the uniform persistence of the complete family of systems over the hull is not to be expected in general, since the appropriate concept of uniform persistence in non-autonomous systems must be a collective one (see~\cite{obsa18} for more details).

Inspired by the treatment given to the scalar model in~\cite{smit}, we now follow his alternative strategy using the exponential ordering  to determine conditions of a different nature on the existence of a unique positive almost periodic solution attracting every other positive solution as $t\to \infty$. When we write the Nicholson systems~\eqref{nicholson delay hull} in the general form $y'(t)=F(\w{\cdot}t,y_t)$, $\w\in\W$, it is easy to  check that $F$ satisfies condition~\ref{F1}. Now, we introduce the cone of positive elements for the exponential ordering. We take an $m\times m$ diagonal matrix $B$ with diagonal entries $-\mu_1, -\mu_2,\ldots, -\mu_m$, for some $\mu_i\geq 0$, $i=1,\ldots,m$. The diagonal structure of $B$ permits to write the positive cone as
\[
K_B= \big\{\varphi\in X \mid \varphi\geq 0  \;\text{and}\; \varphi_i(t)\,e^{\mu_i t}\geq \varphi_i(s)\,e^{\mu_i s},\;-r_i\leq s\leq t\leq 0\,, \;1\leq i\leq m\big\}\,.
\]

To find a sufficient condition relating the coefficients of the Nicholson systems and the matrix $B$ leading to the quasimonotone condition~\ref{F2}, we denote by $D(\w)$  and $\beta(\w)$ the diagonal matrices with diagonal entries $d_1(\w),\ldots,d_m(\w)$ and $\beta_1(\w),\ldots,\beta_m(\w)$, respectively, and by $A(\w)$ the matrix $[a_{ij}(\w)]$. We require that for each $\w\in\W$, whenever $\phi\leq_B \psi$, that is, $\psi-\phi\in K_B$, it holds that
\begin{align*}
F(\w,\psi)-F(\w,\phi)-B\,(\psi(0)-\phi(0))= [-D(\w)+A(\w)-B]\,(\psi(0)-\phi(0))&\\
+\beta(\w)\,\big(\psi_i(-r_i)\,e^{-c_i(\w)\,\psi_i(-r_i)}-\phi_i(-r_i)\,e^{-c_i(\w)\,\phi_i(-r_i)}\big)_{1\leq i\leq m}&\geq 0\,.
\end{align*}
For each component, applying in the first inequality the mean value theorem to the map $h_i(y)=y\,e^{-c_i(\w)\,y}$,  and the fact that $\phi\leq_B \psi$ in the second one, we have that
\begin{align*}
\psi_i(-r_i)\,e^{-c_i(\w)\,\psi_i(-r_i)}-\phi_i(-r_i)\,e^{-c_i(\w)\,\phi_i(-r_i)}&\geq \frac{-1}{e^2}\,(\psi_i(-r_i)-\phi_i(-r_i))\\
&\geq \frac{-1}{e^2}\,e^{\mu_ir_i}(\psi_i(0)-\phi_i(0))\,.
\end{align*}
From here it is easy to get this sufficient condition for~\ref{F2}:
\begin{equation}\label{condicion suf}
\beta_i(\w)\leq (-d_i(\w)+\mu_i)\,e^{2-\mu_ir_i}  \quad \text{for}\;\,\w\in\W \;  \text{ and }\;1\leq i\leq m\,.
\end{equation}
Note that the restriction $d_i(\w)\leq \mu_i$ for all $\w\in\W$ and $1\leq i\leq m$ is implicit in the previous condition and that it is independent of the migration terms $a_{ij}(\w)$ as well as of the coefficients $c_i(\w)$. Next, we determine some precise relations so that~\eqref{condicion suf} holds for some positive $\mu_1,\ldots,\mu_m$. Let us introduce some notation:
\begin{equation}\label{eq:notacion}
d_i^+:=\sup_{t\in\R}\wit d_i(t)\,,\quad  \beta_i^+:=\sup_{t\in\R}\wit \beta_i(t) \,,\quad \text{for}\;\, 1\leq i\leq m\,.
\end{equation}
By the hull construction, $d_i^+=\sup_{\w\in\W}d_i(\w)$ and $\beta_i^+=\sup_{\w\in\W}\beta_i(\w)$ for all $1\leq i\leq m$.

\begin{theorem}\label{teor-nicho}
Assume that the Nicholson system~\eqref{nicholson delay} satisfies conditions~\ref{a1}--\ref{a6} and for each $1\leq i\leq m$ let
$\beta_i^+$ and $d_i^+$ be defined as in~\eqref{eq:notacion}. Let $B$ be the diagonal matrix with diagonal elements $-\mu_1,\ldots,-\mu_m$, for
\begin{equation}\label{mui}
\mu_i= \frac{1}{r_i}\,\ln\left(\frac{e^2}{r_i\,\beta_i^+}\right)\quad \text{for each}\;\, 1\leq i\leq m\,.
\end{equation}
 \begin{itemize}
\item[\rm{(i)}] Whenever
\begin{equation}\label{pequeno}
r_i\,\beta_i^+\,e^{d_i^+ r_i} \leq e  \quad \text{for all}\;\, 1\leq i\leq m\,,
\end{equation}
the matrix $B$ has negative diagonal entries and  condition~\ref{F2} holds, whence the induced skew-product semiflow $\tau$ is monotone for  $\leq_B$.

\item[\rm{(ii)}] Whenever
\begin{equation}\label{pequeno estricto}
r_i\,\beta_i^+\,e^{d_i^+ r_i} < e \quad \text{for all}\;\, 1\leq i\leq m\,,
\end{equation}
condition~\ref{F3} holds and besides, for  $r:=\max(r_1,\ldots,r_m)${\rm :}
\begin{itemize}
\item[\rm(ii.1)] if $\w\in\W$ and $\phi\geq 0$, then $y_t(\w,\phi)\geq_B 0$ for $t\geq r$;
\item[\rm(ii.2)] if $\w\in\W$ and $\phi\geq 0$ with $\phi(0)\gg 0$, then  $y_t(\w,\phi)\gg_B 0$ for $t\geq 2r$.
\end{itemize}
\end{itemize}
\end{theorem}
\begin{proof}
Take into account the sufficient condition~\eqref{condicion suf}  for~\ref{F2} and the fact that, for each $i\in\{1,\ldots,m\}$, by the hull construction,   $-d_i(\w)+\mu_i-e^{\mu_ir_i-2}\beta_i(\w)\geq -d_i^++\mu_i-e^{\mu_ir_i-2} \beta_i^+$ for all $\w\in\W$. Then, it suffices to do an analytical study of the map
\[
f_i(\mu)= - d_i^+ +\mu -\frac{1}{e^2}\,\beta_i^+\,e^{\mu r_i}\,,\quad \mu\geq 0\,,
\]
to see that $f_i(\mu_i)\geq 0$ provided that~\eqref{pequeno} holds. More precisely, $\mu_i > 0$ is the point where the map $f_i$ reaches its maximum value on the positive semiaxis and condition~\eqref{pequeno} guarantees that this maximum value is greater than or equal to $0$.

When~\eqref{pequeno estricto} holds, then the maximum value of $f_i$ on the positive semiaxis is positive, the inequalities in~\eqref{condicion suf} are strict, and thus
condition~\ref{F3} holds.  Now, take  $\w\in\W$ and $\phi\geq 0$.
Then, for $t\geq r$, $y_t(\w,\phi)$ is a smooth map, $y_t(\w,\phi)\geq 0$ and
\[
y_i'(t+s,\w,\phi)\geq -d_i(\w{\cdot}(t+s))\,y_i(t+s,\w,\phi)\geq -\mu_i\, y_i(t+s,\w,\phi)\quad \text{for all } s\in [-r_i,0]
\]
and $1\leq i\leq m$. Hence, $y_t(\w,\phi)\geq_B 0$ for $t\geq r$. Finally, if also $\phi(0)\gg 0$, as we have already mentioned, then $y(t,\w,\phi)\gg 0$ for all $t\geq 0$. To finish, we can argue exactly as in the proof of Proposition~\ref{prop monot fuerte}(i).
\end{proof}
Note that conditions~\eqref{pequeno} and~\eqref{pequeno estricto}  are satisfied whenever the delays are small enough.
Hereafter, when one of these conditions is assumed, the exponential ordering $\leq_B$ for the diagonal matrix $B$ with diagonal entries $-\mu_1, \ldots, -\mu_n$ is considered, in accordance with~\eqref{mui}.

As a consequence of the previous result, we obtain two  corollaries under condition~\eqref{pequeno estricto}. Firstly, if a minimal set lies inside the interior of the standard positive cone, in fact it is inside the smaller set $\Int \wit K_B$ for the exponential ordering. Secondly,  for Nicholson systems the property of uniform persistence of the semiflow with respect to  the exponential order implies the same property for the usual order. Eventually, we will see that both are  equivalent. Recall that these properties of persistence are  in general unrelated.
\begin{coro}\label{coro-1}
Assume that $\eqref{pequeno estricto}$ holds for the Nicholson system $\eqref{nicholson delay}$. If there exists a strongly positive minimal set $M\gg 0$ for $\tau$, then $M\gg_B 0$.
\end{coro}
\begin{proof}
Just note that inside a minimal set there are backward extensions of the semitrajectories. Then, when we move backwards we remain in the interior of the standard positive cone, and when we come back forwards we enter the interior of the exponential ordering cone, by Theorem~\ref{teor-nicho}(ii.2).
\end{proof}
\begin{coro}\label{coro-2}
Assume that $\eqref{pequeno estricto}$ holds for the Nicholson system $\eqref{nicholson delay}$. If the induced semiflow $\tau$ is uniformly persistent in $\Int\wit K_B$, then it is also uniformly persistent in $\Int X_+$.
\end{coro}
\begin{proof}
Let us assume that $\tau$ is uniformly persistent in $\Int\wit K_B$ and let us see that it is also uniformly persistent in $\Int X_+$. Let $\psi\gg_B 0$ be the one in Definition~\ref{defi persistence hull}(ii) and take $\w\in\W$ and $\varphi\gg 0$. By Theorem~\ref{teor-nicho}(ii) we know that $y_{2r}(\w,\varphi)\gg_B 0$, so that there exists a $t_0=t_0(\w,\varphi,2r)$ such that
     $y_t(\w{\cdot}2r, y_{2r}(\w,\varphi))=y_{t+2r}(\w,\varphi)\geq_B \psi$ for all $t\geq t_0$. By the definition of $\leq_B$ in particular $y_{t+2r}(\w,\varphi)\geq \psi\gg 0$, and we are done.
\end{proof}
The next result for the exponential ordering is an application of Theorem~\ref{th:main}. Note that, if we are working with the exponential order $\leq_B$, it seems natural to assume the uniform persistence with respect to its positive cone. However, it is enough to assume this property for the usual order, which, by the previous corollary, is an a priori weaker assumption.
\begin{theorem}\label{teor Main Nicholson}
Assume that $\eqref{pequeno estricto}$ holds for the Nicholson system $\eqref{nicholson delay}$, which satisfies conditions~\ref{a1}--\ref{a6} and is uniformly persistent at $0$. Then, there exists a unique minimal set $M\gg 0$ which is a copy of the base, that is, $M=\{(\w,b(\w))\mid \w\in\ \W\}$ for a continuous map $b:\Omega\to \Int X_+$, and
\begin{equation}\label{main-limit}
\lim_{t\to\infty}\|y_t(\w,\phi)-b(\w{\cdot}t)\|_\infty=0 \quad \text{whenever } \w \in \W\,, \; \phi\geq 0\text{ with }\phi(0)\gg 0\,.
\end{equation}
In other words, for each of the systems in the family $\eqref{nicholson delay hull}$ there is a unique strongly positive almost periodic solution and, whenever at time $0$ all the patches are inhabited,  the population evolution is asymptotically almost periodic.
\end{theorem}
\begin{proof}
In order to apply Theorem~\ref{th:main}, we write the systems~\eqref{nicholson delay hull} in the general form $y'(t)=F(\w{\cdot}t,y_t)$, $\w\in\W$. Then, condition~\ref{F1} holds and, by Theorem~\ref{teor-nicho}(ii), \ref{F3} (and hence~\ref{F2}) follows from~\eqref{pequeno estricto}. It is also easy to see  that the sublinearity condition~\ref{F4} actually holds for all $\psi\geq 0$, so that by Proposition~\ref{prop:sublinear} the semiflow is sublinear for the exponential ordering.

As mentioned before,~\cite[Theorem~3.4]{obsa18} says that Definition~\ref{defi persistence Nicholson} is equivalent to Definition~\ref{defi persistence hull}(i). Since the solutions are ultimately bounded, taking any pair $(\w,\phi)\in\W\times \Int X_+$, its omega-limit set contains a minimal set which  satisfies $M\gg 0$ by the uniform persistence property. But then, by Corollary~\ref{coro-1}, $M\gg_B 0$ and all the semitrajectories inside $M$ satisfy condition~\ref{F5}.

Finally, for condition~\ref{F6}, let $\w_1$ be the map in $\W$ determining the initial Nicholson system~\eqref{nicholson delay}, which we write as
$y_i'(t)=f_i(t,y_t)$, $1\leq i\leq m$, for the maps $f_i:\R\times X\to \R$ given in \eqref{eq:f}.
By assumptions~\ref{a4}--\ref{a5}, for each fixed $t$, the real map $h(x)=\wit \beta_{i}(t)\,xe^{-\wit c_i(t)\,x}$ is strictly sublinear for $x>0$. Thus, it is easy to check that, for each $\lambda\in (0,1)$ and $t\geq 0$,
\begin{equation}\label{F st sublinear}
\phi\geq 0\; \text{with}\;\,\phi_i(-r_i)>0\,,\; 1\leq i\leq m \;\Rightarrow\; \lambda\, f(t,\phi)\ll f(t,\lambda\,\phi)\,.
\end{equation}
We claim that~\ref{F6} holds for $\w_1$ and every $t_1>r:=\max(r_1,\ldots,r_m)$. To check it, take $\psi\gg_B 0$ and $\lambda \in(0,1)$.
Having in mind the expression of  $\Int\wit K_B$, let us first  check that $\lambda\, y_t(\w_1,\psi)\ll y_t(\w_1,\lambda\, \psi)$ for all $t>r$ by comparing the solutions.
It suffices to see that $\lambda\, y(t,\w_1,\psi)\ll y(t,\w_1,\lambda\, \psi)$ for all $t>0$. Consider $\wit r:=\min(r_1,\ldots,r_m)>0$ and let $z(t):=\lambda\, y(t,\w_1,\psi)$, $t\geq 0$. Then, since in particular $\psi\gg 0$, again by strict sublinearity,  for $t\in (0,\wit r]$ and $1\leq i\leq m$ we have that
\begin{align*}
z_i'(t)&=-\wit d_i(t)\,z_i(t) +\des \sum_{j=1}^m \wit a_{ij}(t)\,z_j(t) + \wit\beta_{i}(t)\,\lambda\,\psi_i(t-r_i)\,e^{-\wit c_i(t)\,\psi_i(t-r_i)}\\
&< -\wit d_i(t)\,z_i(t) +\des \sum_{j=1}^m \wit a_{ij}(t)\,z_j(t) + \wit\beta_{i}(t)\,\lambda\,\psi_i(t-r_i)\,e^{-\lambda\,\wit c_i(t)\,\psi_i(t-r_i)}\,.
\end{align*}
A standard comparison theorem for cooperative  ODEs (thanks to~\ref{a3}) implies that $z(t)\ll y(t,\w_1,\lambda\, \psi)$ for $t\in (0,\wit r]$. Iterating the procedure interval by interval, we can assert that $\lambda\, y(t,\w_1,\psi)\ll y(t,\w_1,\lambda\, \psi)$ for all $t>0$, as we wanted.

Secondly, note that, since  $\psi\gg 0$, we know that  $y(t,\w_1,\psi)\gg 0$ for all $t\geq 0$, so that $y_t(\w_1,\psi)\gg 0$ for all $t\geq 0$. Then,  for all $t>0$ and $1\leq i\leq m$ we write
\begin{align*}
 &y_i'(t,\w_1,\lambda\,\psi) -\lambda \,y_i'(t,\w_1,\psi)+\mu_i(y_i(t,\w_1,\lambda\,\psi) -\lambda\, y_i(t,\w_1,\psi)) \\
 &\qquad = f_i(t,y_t(\w_1,\lambda\,\psi))- \lambda\, f_i(t,y_t(\w_1,\psi))+\mu_i(y_i(t,\w_1,\lambda\,\psi) -\lambda\, y_i(t,\w_1,\psi))\\
 &\qquad >  f_i(t,y_t(\w_1,\lambda\,\psi))-  f_i(t,\lambda\, y_t(\w_1,\psi))+\mu_i(y_i(t,\w_1,\lambda\,\psi) -\lambda\, y_i(t,\w_1,\psi))
 \geq 0\,,
\end{align*}
where~\eqref{F st sublinear} applied to $y_t(\w_1,\psi)\gg 0$  has been used in the first inequality, and the sublinearity of the semiflow for $\leq_B$ permits to apply~\ref{F2} to justify  the second inequality. Then, by a continuity argument we deduce that $y_t(\w_1,\lambda\,\psi)\gg_B \lambda \,y_t(\w_1,\psi)$ for all $t>r$, whence~\ref{F6} holds for $\w_1$ and every  $t_1> r$, as claimed.

Therefore, we can apply  Theorem~\ref{th:main} saying that there exists a unique strongly positive minimal set $M\gg_B 0$ which is a copy of the base and attracts all the semiorbits starting in $\W\times\Int \wit K_B$ as $t\to\infty$ for the $\|\cdot\|_L$ norm. By Corollary~\ref{coro-1}, $M$ is also the unique minimal set in $\Int X_+$. Finally,~\eqref{main-limit} follows from Theorem~\ref{teor-nicho}(ii.2), since the orbit of each  $(\w,\phi)$ with  $\phi\geq 0$ and $\phi(0)\gg 0$ enters $\W\times\Int \wit K_B$ after time $2r$ and it is then attracted to the minimal set  forwards in time. Since  the immersion $X_L\hookrightarrow X$ is continuous,  the proof is finished.
\end{proof}

\begin{coro}\label{coro-3}
Assume that $\eqref{pequeno estricto}$ holds for the Nicholson system $\eqref{nicholson delay}$. The induced semiflow $\tau$ is uniformly persistent in $\Int\wit K_B$ if and only if it is uniformly persistent in $\Int X_+$.
\end{coro}
\begin{proof}
One implication is Corollary~\ref{coro-2}. Assume that $\tau$ is uniformly persistent in $\Int X_+$. Thanks to condition~\eqref{pequeno estricto}, Theorem~\ref{teor Main Nicholson} applies so that the set $M=\{(\w,b(\w))\mid \w\in\W\}\gg_B 0$ is the attractor of  positive solutions. In particular, there exists a $\psi\gg_B 0$ such that $b(\w)\geq_B \psi$ for all $\in\W$. Then, taking any $\w\in\W$ and $\varphi\gg_B 0$,  the semiorbit $\{(\w{\cdot}t,y_t(\w,\varphi))\mid t\geq 0\}$ approaches $M$ as $t\to\infty$ and therefore, there exists a $t_0>0$ such that $y_t(\w,\varphi)\geq_B \psi/2$ for all $t\geq t_0$. The proof is finished.
\end{proof}
If some of the maps $\wit d_i(t)$ have a big oscillation with respect to their mean values,
\[
d_{i0} := \lim_{t\to\infty} \frac{1}{t}\int_0^t \wit d_i(s)\,ds\quad \text{for}\;\, 1\leq i\leq m\,,
\]
condition~\eqref{pequeno estricto},  $r_i\,\beta_i^+\,e^{d_i^+ r_i} < e$, might force  too strong a small delay condition. In the next result we relax condition~\eqref{pequeno estricto}.
\begin{theorem}\label{teor cv d}
Assume that the Nicholson system $\eqref{nicholson delay}$ satisfies conditions~\ref{a1}--\ref{a6} and it is uniformly persistent at $0$. If  it also holds that
\begin{equation}\label{eq:condicion cv d}
r_i\,\sup\left\{\wit\beta_i(t)\,e^{\int_{t-r_i}^t \wit d_i(s)\,ds} \,\big|\, t\in\R\right\} < e\quad \text{for each}\;\, 1\leq i\leq m\,,
\end{equation}
then, there exists a unique positive almost periodic solution of~\eqref{nicholson delay} which asymptotically attracts every other positive solution; more precisely, it attracts  every other solution $y(t,\varphi)$ with initial value $\varphi\geq 0$ such that $\varphi(0)\gg 0$.
\end{theorem}
\begin{proof}
Consider the family of systems~\eqref{nicholson delay hull} over the hull $\W$, which is minimal and uniquely ergodic because of the almost periodicity assumption~\ref{a1}. Let $\nu$ be the unique ergodic measure on $\W$.
Roughly speaking, the idea is to make a joint change of variables which for each $\w\in \W$ takes the almost periodic coefficients $d_i(\w{\cdot}t)$, $1\leq i\leq m$ into their respective mean values $d_{i0}$, or at least arbitrarily close to them. Recall that, by Birkhoff's ergodic theorem, $d_{i0}=\int_\W d_i\,d\nu$ for $1\leq i\leq m$.  Consider the Banach space $C_0(\W)=\big\{a\in C(\W)\mid \int_\W a\,d\nu=0\big\}$ of continuous maps on $\W$ with null mean value, and its  vector subspace $BP(\W)=\{ a\in C_0(\W) \mid a \,\text{ has a bounded  primitive}\}$.  In the case of periodic coefficients, $C_0(\W)=BP(\W)$, whereas, if $\W$ is aperiodic, then $BP(\W)\subsetneq C_0(\W)$ and it is a dense subset of first category. See Campos et al.~\cite[Lemma~5.1]{caot} for this result in a more general setting, but note that the almost periodic case was already proved by Johnson~\cite{john80}.

We could skip the first case we are going to consider, but it helps to understand the condition required in~\eqref{eq:condicion cv d}. Let us first assume that the maps $\widehat d_i (\w) = d_i(\w)-d_{i0}$, $\w\in\W$ for $1\leq i\leq m$, which are in $C_0(\W)$, admit bounded primitives, that is, there exist maps $h_i\in C(\W)$,  $1\leq i\leq m$, such that
\begin{equation}\label{eq:primitiva}
\int_0^t \widehat d_i (\w{\cdot}s)\,ds = h_i(\w{\cdot}t)-h_i(\w) \quad\text{for all}\;\,\w\in\W\,,\; t\in \R\,.
\end{equation}
Then, in \eqref{nicholson delay hull}, for each  $\w\in\W$  we make the change of variables $z_i(t)=e^{h_i(\w{\cdot}t)}\,y_i(t)$, $1\leq i\leq m$, which preserves the boundedness of the solutions and the property of uniform persistence. The transformed systems have a similar structure to the Nicholson systems, namely, for each $\w\in\W$,
\begin{equation*}
z'_i(t)=- d_{i0}\,z_i(t) +\des \sum_{j=1}^m  a_{ij}^*(\w{\cdot}t)\,z_j(t) + \beta_{i}^*(\w{\cdot}t)\,z_i(t-r_i)\,e^{-c_i^*(\w{\cdot}t)\,z_i(t-r_i)}\,
\end{equation*}
for the new coefficients
\begin{align*}
a_{ij}^*(\w):&= a_{ij}(\w)\,e^{h_i(\w)-h_j(\w)}\,, \\
\beta_{i}^*(\w):&= \beta_{i}(\w)\,e^{h_i(\w)-h_i(\w{\cdot}(-r_i))}\,,\\
c_{i}^*(\w):&= c_{i}(\w)\,e^{-h_i(\w{\cdot}(-r_i))}\,,
\end{align*}
and it is easy to check that the transformed system of~\eqref{nicholson delay} also satisfies~\ref{a1}--\ref{a5}, although it might not satisfy assumption~\ref{a6}. However, this is unimportant because, even if this condition has a biological meaning in the original system,  its analytical implications --it is basically used to prove the ultimate  boundedness of solutions-- still hold, by the expression of the change of variables,  the fact that $\W$ is compact, and $h_i$ are continuous maps.

Therefore, we can apply Theorem~\ref{teor Main Nicholson} to the new systems, provided that the monotonicity condition for the exponential ordering~\eqref{pequeno estricto} holds, but that is exactly condition~\eqref{eq:condicion cv d} in the hypotheses. To check it, let $\w_1\in\W$ be the one determining the initial Nicholson system~\eqref{nicholson delay} and note that now the parameters involved for each $1\leq i\leq m$ are  $ (d_i^*)^+= d_{i0}$ and
\[
(\beta_i^*)^+ := \sup\big\{ \beta_i^*(\w_1{\cdot}t) \mid t\in\R\big\} =\sup\left\{  \beta_{i}(\w_1{\cdot}t)\,e^{h_i(\w_1{\cdot}t)-h_i(\w_1{\cdot}(t-r_i))}\,\big|\, t\in\R\right\}\,,
\]
which can be rewritten as $(\beta_i^*)^+ =\sup\left\{\wit\beta_i(t)\,e^{\int_{t-r_i}^t \wit d_i(s)\,ds -d_{i0}r_i} \,\big|\, t\in\R\right\}$, by~\eqref{eq:primitiva}. As a consequence, for the transformed systems there exists a unique minimal set $M^*=\{(\w,b^*(\w))\mid \w\in\W\}$ for a continuous map $b^*:\W\to \Int X_+$, which is an attractor of the semitrajectory of every pair $(\w,\phi)$ provided that $\phi\geq 0$ and $\phi(0)\gg 0$. Reversing the change of variables, we see that the original systems~\eqref{nicholson delay hull} also have a  unique minimal set $M=\{(\w,b(\w))\mid \w\in\W\}$ for a continuous map $b:\W\to \Int X_+$  which is also an attractor of the semitrajectory of every pair $(\w,\phi)$ provided that $\phi\geq 0$ and $\phi(0)\gg 0$.

To finish the proof it remains to assume that some of the maps $\widehat d_i (\w) = d_i(\w)-d_{i0}$, $\w\in\W$, for $1\leq i\leq m$, which are in $C_0(\W)$, do not admit a bounded primitive. Then, since the subspace $BP(\W)$ is dense in $C_0(\W)$ and~\eqref{eq:condicion cv d} holds, given any $\ep>0$ with the following restrictions:
\begin{itemize}
\item[($\ep$1)] $\ep < d_{i0}/2\,$,
\item[($\ep$2)] $r_i\,\sup\left\{\wit\beta_i(t)\,e^{\int_{t-r_i}^t \wit d_i(s)\,ds} \,\big|\, t\in\R\right\} < e^{1-2\ep r_i}\,$ for $\,1\leq i\leq m$\,,
\end{itemize}
there exist maps $\widehat \delta_i\in BP(\W)$, $\,1\leq i\leq m$ (which depend on $\ep$), such that $\sup\big\{|\widehat d_i(\w)-\widehat \delta_i(\w)|\mid \w\in\W\big\}\leq \ep$. This time let $h_i\in C(\W)$ be a primitive of $\widehat \delta_i$, that is,
\begin{equation}\label{eq:primitiva 2}
\int_0^t \widehat \delta_i (\w{\cdot}s)\,ds = h_i(\w{\cdot}t)-h_i(\w)  \quad\text{for all}\;\,\w\in\W\,,\; t\in \R\,.
\end{equation}
The same change of variables as before, $z_i(t)=e^{h_i(\w{\cdot}t)}\,y_i(t)$, $1\leq i\leq m$, transforms the systems into the new family
\begin{multline*}
z'_i(t)=-\big[d_{i0}+(\widehat d_i(\w{\cdot}t)-\widehat \delta_i(\w{\cdot}t))\big]\,z_i(t) \\ +\des \sum_{j=1}^m  a_{ij}^*(\w{\cdot}t)\,z_j(t) + \beta_{i}^*(\w{\cdot}t)\,z_i(t-r_i)\,e^{-c_i^*(\w{\cdot}t)\,z_i(t-r_i)}\,,\quad \w\in\W\,,
\end{multline*}
where the coefficients $a_{ij}^*(\w)$, $\beta_{i}^*(\w)$ and $c_{i}^*(\w)$ have the same expressions as before. It is easy to see that under restriction ($\ep$1) the new transformed system of~\eqref{nicholson delay} also satisfies~\ref{a1}--\ref{a5}. Besides, taking into account~\eqref{eq:primitiva 2}  and doing similar calculations to the ones in the previous situation, we see that,  thanks to restriction   ($\ep$2), relation~\eqref{pequeno estricto} also holds for this transformed system. Hence, Theorem~\ref{teor Main Nicholson} applies and the proof is finished as before, just by reversing the change of variables.
\end{proof}
\begin{notas}
\par
\rm{\bf{1.}} For each $i=1,\ldots,m$, note that, if the map $\wit d_i(t)$ is periodic and the delay $r_i$ is a multiple of the period, then $\int_{t-r_i}^t \wit d_i(s)\,ds= d_{i0}r_i$ for all $t\in \R$, and  condition~\eqref{eq:condicion cv d} then reads $r_i\,\beta_i^+\,e^{ d_{i0}r_i}  < e$. This last condition allows for bigger delays (and bigger periods too) rather than $r_i\,\beta_i^+\,e^{ d_i^+r_i} < e$ in~\eqref{pequeno estricto}, specially for periodic maps with a big oscillation. Analogously, if the map $\wit d_i(t)$ is almost periodic, given an $\varepsilon>0$ as small as wanted,  there is a relatively dense set in $(0,\infty)$ of delays $r_i$ such that $\big|\int_{t-r_i}^t \wit d_i(s)\,ds -  d_{i0}r_i\big| <\varepsilon$ for all $t\in \R$, so that  $r_i\,\beta_i^+\,e^{ d_{i0}r_i+\varepsilon}  < e$ implies  condition~\eqref{eq:condicion cv d}.  Whereas the statement for the periodic case is immediate, in the almost periodic case we have to argue as in the proof of the last result, taking into account that continuous maps on $\W$ with null mean value can be approximated by maps with a continuous primitive, and also the fact that an almost periodic map has a relatively dense set of $\delta$-periods for each $\delta>0$.

\rm{\bf{2.}} The same technique used in the proof of the previous theorem allows to improve the small delay conditions imposed in the literature for the existence of the so-called {\em special solutions\/} of FDEs, when applied to the scalar Nicholson  equation $x'(t)=-\wit d(t)\,x(t)+\wit \beta(t)\,x(t-r)\,e^{-\wit c(t)\,x(t-r)}$ with almost periodic coefficients $\wit d(t)$, $\wit \beta(t)$ and $\wit c(t)$.  Special solutions are globally defined solutions (defined on the whole line $\R$) of delay equations, which are solutions of associated ODEs.  This theory originated in the 1960s in some works by Ryabov~\cite{ryab}. See also~\cite{driv},~\cite{pitu}, and the references therein. More precisely, fix a general scalar delay equation $x'(t)=f(t,x_t)$ with delay $r$, such that $f$ is continuous, satisfies $\sup_{t\leq 0}|f(t,0)|\, e^{t/r} < \infty$ and is globally Lipschitz, that is, $|f(t,\phi)-f(t,\psi)|\leq L\,\|\phi-\psi\|$ for all $t\in\R$, $\phi,\psi\in C([-r,0])$. Then, provided that the delay satisfies $L\,r\,e <1$, for each $t_0, x_0\in \R$, there exists a unique solution $x(t)$ of $x'=f(t,x_t)$, $t\in \R$,  such that $x(t_0)=x_0$ and  $\sup_{t\leq 0}|x(t)|\, e^{t/r} < \infty$. This solution is the special solution. In particular the solutions in the global attractor, when it exists, are globally defined and bounded, and thus are special solutions.

For the Nicholson equation, $L=d^+ + \beta^+$ and the previous small delay condition reads $(d^+ + \beta^+)\,r\,e<1$.
By looking at the family of equations over the hull and performing the change of variables $z(t)=e^{h(\w{\cdot}t)}\,x(t)$ for an appropriate $h\in C(\W)$, which clearly preserves special solutions,  we can now improve the small delay condition for the existence of special solutions:
\[
\left(d_0+\sup\left\{\wit\beta(t)\,e^{\int_{t-r}^t \wit d(s)\,ds -d_{0}r} \,\big|\, t\in\R\right\}\right)\,r\,e <1\,,
\]
where $d_0$ is the mean value of $\wit d(t)$.   If in particular  $\wit d(t)$ is periodic and the delay $r$ is a multiple of the period, then this condition is just $(d_0+ \beta^+)\,r\,e<1$.

The same happens when we look at more recent results by Pituk~\cite{pitu} for scalar autonomous delay equations. By introducing the exponential ordering $\leq_{\mu}$ for a $\mu>0$, \cite[Theorem~5.2]{pitu} offers new conditions  to guarantee the existence of special solutions for $x'(t)=f(x_t)$, this time with the growth condition $\sup_{t\leq 0}|x(t)|\, e^{\mu t} < \infty$. It is easy to check that~\cite[Theorem~5.2]{pitu} is also true  for time-dependent equations $x'(t)=f(t,x_t)$ if condition (5.1) therein is independent of $t$ and besides $\sup_{t\leq 0}|f(t,0)|\, e^{\mu t} < \infty$. The condition to apply this result to our Nicholson equation is $r\,\beta^+  e^{d^+ r}<e^{-1}$.
Once more performing the previous change of variables, we can improve this condition into
\[
r\,\sup\left\{\wit\beta(t)\,e^{\int_{t-r}^t \wit d(s)\,ds} \,\big|\, t\in\R\right\}<\frac{1}{e}\,,
\]
which in the aforementioned case of $\wit d(t)$ periodic reduces to $r\,\beta^+ e^{d_0 r}<e^{-1}$.
\end{notas}

\end{document}